\documentclass[12pt,letterpaper]{article}
\usepackage[utf8]{inputenc}

\usepackage{amsfonts,amssymb,amsmath,amsthm} 
\usepackage{graphicx} 
\usepackage{float}
\usepackage[colorlinks=true,linkcolor=blue,citecolor=green]{hyperref}

\usepackage{xcolor}
\definecolor{darkgreen}{rgb}{0,0.4,0}
\usepackage{cite} %<-les cites sont automatiquement dans l'ordre de la biblio [13,26,45] et non pas [23,13,45], voir si on veut cette option au final ou pas...
\colorlet{ligthgray}{gray!40}
\usepackage{euscript} %for the \F that Rafik likes

\newtheorem{Th}{Theorem}
\newtheorem{Lem}{Lemma}

\newcommand{\Z}{\mathbb{Z}}
\newcommand{\R}{\mathbb{R}}

\newcommand{\F}{\EuScript{F}}

\renewcommand{\P}{\mathbb{P}}
\newcommand{\E}[1]{{\mathbb E}\!\left[#1\right]}

\def\build#1_#2^#3{\mathrel{\mathop{#1}\limits_{#2}^{#3}}}
\def\converge#1#2#3#4{\build\hbox to #1mm{\rightarrowfill}_{#2\rightarrow #3}^{\hbox{\scriptsize #4}}}
\font\dsrom=dsrom10 scaled 1200
\def \ind{\textrm{\dsrom{1}}}
\numberwithin{equation}{section}

\begin{document}

\title{On a class of unbalanced step-reinforced random walks}
\author{R. Aguech$^{1}$, S. Ben Hariz$^{2}$, M. El Machkouri$^{3}$, Y. Faouzi$^{4}$ \\[3mm]
	$^{1}$ Department of Statistics and Operations Research, \\ King Saud University, Riyadh, Saudi Arabia.\\ [3mm]
	$^{2}$ Department of Mathematics, Le Mans University, France.\\[3mm]
	$^{3,4}$ EMINES - School of Industrial Management, \\
	UM6P - University Mohammed VI Polytechnic, \\
	Benguerir, Morocco.}	
	\vspace{-1mm}	
	\maketitle
	
	\begin{abstract}
A step-reinforced random walk is a discrete-time stochastic process with long-range dependence. At each step, with a fixed probability $\alpha$, the so-called positively step-reinforced random walk repeats one of its previous steps, chosen randomly and uniformly from its entire history. Alternatively, with probability $1-\alpha$, it makes an independent move. For the so-called negatively step-reinforced random walk, the process is similar, but any repeated step is taken with its direction reversed. These random walks have been introduced respectively by Simon (1955) and Bertoin (2024) and are sometimes refered to the self-confident step-reinforced random walk and the counterbalanced step-reinforced random walk respectively. In this work, we introduce a new class of unbalanced step-reinforced random walks for which we prove the strong law of large numbers and the central limit theorem. In particular, our work provides a unified treatment of the elephant random walk introduced by Schutz and Trimper (2004) and the positively and negatively step-reinforced random walks.\\
\\
\vspace{+2mm}
\noindent \textit{AMS Subject Classifications (2020)}: 60G42, 60F05, 60F25. .\\
\noindent \textit{Keywords}: elephant random walk, step-reinforced random walk, central limit theorem, asymptotic normality, martingale, quasi-martingale, law of large numbers.
	\end{abstract}
	
\section{Introduction and main results}
The most famous example of step-reinforced random walk is the so-called elephant random walk (ERW) on $\Z$ introduced by Schütz and Trimper \cite{Schutz--Trimper2004}. Its name is motivated by the well-known idea that elephants have excellent memories and are able to remember at any time all the places they have visited. At time \( n = 0 \), the elephant is at position $Z_0 = 0$. At time \( n = 1 \), the elephant moves toward $1$ (i.e. ``one step to the right'') with probability $s$ and toward $-1$ (i.e. ``one step to the left'') with probability $1-s$ where $s\in[ 0, 1]$ is a fixed parameter. Thus, the position of the elephant at time \( n = 1 \) is given by $Z_1$ where $Z_1$ is a Rademacher random variable satisfying $\P(Z_1=1)=1-\P(Z_1=-1)=s$. Let  $n\geqslant 1$ be a fixed integer and $U_n$ be an integer chosen uniformly at random from the set $\{1, 2, \ldots, n\}$. Then, the $(n+1)$th step $Z_{n+1}$ of the elephant is defined by
	$$
	Z_{n+1} =\left\{
	\begin{array}{ccc}
		Z_{U_n} & \text{with probability} &  p \\
		-Z_{U_n} & \text{with probability} & 1 - p
	\end{array}\right.
	$$
	where $p\in [0,1]$ is the memory parameter of the ERW model. Thus, for $n \geq 2$, the position of the elephant at time \( n \) is
	$\sum_{i=1}^n Z_i$. The description of the asymptotic behavior of the ERW has motivated a great deal of work in the last two decades. In particular, it has been shown (see \cite{Schutz--Trimper2004}) that the dynamics of the ERW is a function of the value of the memory parameter $p$. More precisely, the ERW exhibits three different regimes called {\em diffusive}, {\em critical} and {\em superdiffusive} depending on whether $p<3/4$, $p=3/4$ or $p>3/4$ respectively. It is not possible to give an exhaustive list of the results obtained to date on the ERW. However, one can recall that by establishing an elegant connection of the ERW model with Pólya-type urns, Baur and Bertoin \cite{Baur--Bertoin2016} obtained an invariance principle (functional central limit theorem) and Bercu \cite{Bercu2018} proposed a powerful approach based on martingale theory allowing to establish many limit theorems for the ERW. 
	A step-reinforced random walk is an extension of the ERW by allowing the law of its first step $Z_1$ to be no longer the Rademacher law. More precisely, let $(\xi_n)_{n\geqslant 1}$ be a sequence of i.i.d. real random variables and $\alpha\in [0,1]$, a positively step-reinforced sequence $(Y_n)_{n\geqslant 1}$ is defined in the following way: $Y_1=\xi_1$ a.s. and for any integer $n\geqslant 2$, with probability $\alpha$, the $n$th step $Y_n$ equals one of the previous steps $Y_1$, $Y_2$,...,$Y_{n-1}$ chosen uniformly whereas $Y_n$ is defined as an independent step $\xi_{n}$  with probability $1-\alpha$. Similarly, the sequence $(Y_n)_{n\geqslant 1}$ is called a negatively step-reinforced sequence if $Y_1=\xi_1$ a.s. and for any integer $n\geqslant 2$, with probability $\alpha$, the $n$th step $Y_{n}$ equals the opposite of one of the previous steps $Y_1$, $Y_2$,...,$Y_{n-1}$ chosen uniformly whereas $Y_n$ equals an independent step $\xi_{n}$ with probability $1-\alpha$. The concept of positively random walks goes back to a basic linear reinforcement algorithm which was introduced a long time ago by Simon \cite{Simon1955}. However, the negatively step-reinforced random walk was introduced recently by Bertoin \cite{Bertoin2024} for which he investigated the law of large numbers and the central limit theorem (see also \cite{Hu2025}, \cite{Hu--Zhang2024} and \cite{Kiss--Veto2022}). In this work, we introduce a new class of step-reinforced random walks in order to provide an unified analysis of the ERW and the positively and negatively step-reinforced random walks. More precisely, if $(p,\alpha)\in[0,1]^2$ is fixed, then the unbalanced step-reinforced random walk $(X_n)_{n\geqslant 1}$ with parameter $(p,\alpha)$ is defined by $X_1=\xi_1$ a.s. and for any integer $n\geqslant 1$,
	\begin{equation}\label{definition_unbalanced_step-reinforced_random_walk}
	X_{n+1}= \left\{
	\begin{array}{ccc}
	X_{U_n} & \text { with probability } &  p \alpha \\ 
	-X_{U_n} & \text { with probability } & (1-p) \alpha \\ 
	\xi_{n+1} & \text { with probability } & 1-\alpha
	\end{array}\right.
	\end{equation}
	where $(U_n)_{n\geqslant 1}$ is a sequence of i.i.d. random variables uniformly distributed on $\{1,...,n\}$.
 One can notice that $(X_n)_{n\geqslant 1}$ reduces to a positively or negatively step-reinforced random walk as soon as $p=1$ or  $p=0$ respectively. Moreover, if $\alpha=1$ and $\P(\xi_1=1)=1-\P(\xi_1=-1)=s$ with $s\in[0,1]$, then $(X_n)_{n\geqslant 1}$ reduces to the ERW introduced in \cite{Schutz--Trimper2004}. In this paper, we investigate the validity of the strong law of large numbers (see Theorems \ref{Loi_forte_des_grands_nombres}, \ref{Loi_des_grands_nombres_Marcinkiewicz-Zygmund} and \ref{Loi_des_grands_nombres_sous_moment_ordre_deux} below) and the asymptotic normality of the unbalanced step-reinforced random walk defined by \eqref{definition_unbalanced_step-reinforced_random_walk} in the standard setting (see Theorem \ref{TLC}) and in the triangular array setting (see Theorem \ref{TLC_bis}). In the sequel, we denote
 \begin{equation}\label{definition_sigma_carre_et_a}
 \sigma^2:=\mu_2-\frac{(1-\alpha)^2\mu_1^2}{(1-a)^2}
 \end{equation}
where $a:=(2p-1)\alpha$ and $\mu_{\ell}:=\mathbb E[\xi_1^{\ell}]$ for any $\ell\in\{1,2\}$. Let $T_n:=\sum_{k=1}^nX_k$ be the position of the unbalanced step-reinforced walker at time $n\geqslant 1$. In the following, we are going to consider only the case $-1\leqslant a<1$ since $T_n=n\xi_1$ a.s. when $a=1$ (i.e. $p=\alpha=1$). We focus also only on the case $\alpha\neq 0$ since the case $\alpha=0$ reduces to the full iid case. Our first result is the following strong law of large numbers.
\begin{Th}\label{Loi_forte_des_grands_nombres}
	If $-1\leqslant a<1$ and $\mathbb E[\vert\xi_1\vert]<+\infty$, then $\frac{T_n}{n}\converge{12}{n}{+\infty}{a.s.}\frac{(1-\alpha) \mu_1}{1-a}$.
\end{Th}
\noindent \textbf{Remark 1}. Theorem \ref{Loi_forte_des_grands_nombres} improves Proposition 1.1 in \cite{Bertoin2024} and Theorem 1.1 in \cite{Hu--Zhang2024} where only positively or negatively step-reinforced random walks were considered. Moreover, their proofs rely on a combinatorial analysis approach using random recursive trees which seems to be no longer suited for unbalanced random walks given by \eqref{definition_unbalanced_step-reinforced_random_walk}. Therefore, we adopt an approach based on martingale theory and truncation arguments (see proof of Theorem \ref{Loi_forte_des_grands_nombres} below).\\
\\
\noindent Assuming a stronger moment assumption, one can derive a rate of convergence in the strong law of large numbers given by the following Marcinkiewicz-Zygmund's type law of large numbers. 
\begin{Th}\label{Loi_des_grands_nombres_Marcinkiewicz-Zygmund}
 Let $(X_{n})_{n\geqslant 1}$ be the unbalanced step-reinforced random walk defined by \eqref{definition_unbalanced_step-reinforced_random_walk} such that $\mathbb{E}[\left\vert \xi_{1}\right\vert^{r}]<+\infty$ for some $1<r<2$.
 	\begin{itemize}
 		\item[i)] If $-1\leqslant a<1/2$, then $n^{1-\frac{1}{r}}\left( \frac{T_{n}}{n}-\frac{(1-\alpha )\mu _{1}}{1-a}\right)\converge{12}{n}{+\infty}{a.s.}0$.
 		\item[ii)] If $a=1/2$, then 
 		$
 		\frac{n^{1-\frac{1}{r}}}{\sqrt{\log n}}\left( \frac{T_{n}}{n}-2(1-\alpha )\mu _{1}\right)\converge{12}{n}{+\infty}{a.s.}0.
 		$
 		\item[iii)] If $1/2<a<1$, then 
 		$
        n^{\frac{3}{2}-\frac{1}{r}-a}\left( \frac{T_{n}}{n}-\frac{(1-\alpha )\mu _{1}}{1-a}\right)\converge{12}{n}{+\infty}{a.s.}0.
 		$
 	\end{itemize}
 \end{Th}
\noindent If $\xi_1$ is square-integrable, then the following version of the law of large numbers holds.
\begin{Th}\label{Loi_des_grands_nombres_sous_moment_ordre_deux}
 Let $(X_{n})_{n\geqslant 1}$ be the unbalanced step-reinforced random walk defined by \eqref{definition_unbalanced_step-reinforced_random_walk} such that $\mathbb{E}[\xi_{1}^2]<+\infty$.
	\begin{itemize}
		\item[i)] If $-1\leqslant a<1/2$, then
		$\frac{\sqrt{n}}{\left(\log n\right)^{\frac{1}{2}+\gamma}}\left(\frac{T_n}{n}-\frac{(1-\alpha)\mu _{1}}{1-a}\right)\converge{12}{n}{+\infty}{a.s.}0
		$ for any $\gamma>0$.
		\item[ii)] If $a=1/2$, then
		$
		\frac{\sqrt{n}}{\sqrt{\log n}\left(\log \log n\right)^{\frac{1}{2}+\gamma}}\left(\frac{T_n}{n}-2(1-\alpha)\mu _{1}\right)\converge{12}{n}{+\infty}{a.s.}0$ for any $\gamma>0$.
		\item[iii)] If $1/2<a<1$, then 
		$
		n^{1-a}\left(\frac{T_n}{n}-\frac{(1-\alpha)\mu _{1}}{1-a}\right)\converge{12}{n}{+\infty}{a.s.}L
		$ 
		where $L$ is a square-integrable random variable whose first two moments are given by
\begin{align*}
&\mathbb{E}[L]=\frac{(\alpha-a) \mu_1}{(1-a) \Gamma(a+1)}\quad\textrm{and}\\
&\mathbb{E}[L^2]=\frac{\mu_2}{\Gamma(2 a+1)}-\frac{(1-\alpha)(1+\alpha-2 a) \mu_1^2}{(1-a)^2\Gamma(2a+1)}+\frac{\sigma^2}{(2 a-1) \Gamma(2 a+1)}-\frac{2(1-\alpha)(\alpha-a) \mu_1^2}{(1-a)^2 \Gamma(2 a+1)} 
\end{align*}
where $\sigma^2$ is given by \eqref{definition_sigma_carre_et_a}. Moreover, $L$ is a degenerate random variable if and only if $\mu_2=\mu_1^2$ and either $\mu_1=0$ or $p=1$.
	\end{itemize}
\end{Th}
\noindent \textbf{Remark 2.} The formulas obtained for the first two moments of $L$ in Theorem \ref{Loi_des_grands_nombres_sous_moment_ordre_deux} generalize the corresponding formulas derived by Bercu \cite[Theorem $3.8$]{Bercu2018} and by Kiss and Vet\H{o} \cite[Theorem $1.2$]{Kiss--Veto2022} in the particular cases of the ERW and the positively step-reinforced random walk, respectively. We conjecture that the random variable $L$ in Theorem \ref{Loi_des_grands_nombres_sous_moment_ordre_deux} is non-Gaussian if $\mu_2>\mu_1^2$. To prove this conjecture, one could follow Bercu's approach in \cite[Remark $3.2$]{Bercu2018} by carrying out the lengthy and technical computations of the skewness and kurtosis of $L$ in the more general case of the unbalanced step-reinforced random walk. This interesting question is deferred to future work.
\\
\\
\noindent Our last two results provides the central limit theorem of the unbalanced step-reinforced random walk in both the standard memory setting and the gradually increasing memory setting. 
 \begin{Th}\label{TLC}
 Let $(X_{n})_{n\geqslant 1}$ be the unbalanced step-reinforced random walk defined by \eqref{definition_unbalanced_step-reinforced_random_walk} such that $\mathbb{E}[ \xi_{1}^2]<+\infty$.
\begin{itemize}
	\item[i)] If $-1\leqslant a<1/2$, then $\sqrt{n}\left(\frac{T_n}{n}-\frac{(1-\alpha)\mu _{1}}{1-a}\right)\converge{12}{n}{+\infty}{\textrm{Law}}\mathcal N\left(0,\frac{\sigma^2}{1-2a}\right).
	$
	\item[ii)] If $a=1/2$, then 
	$
	\frac{\sqrt{n}}{\sqrt{\log n}}\left(\frac{T_n}{n}-2(1-\alpha)\mu _{1}\right)\converge{12}{n}{+\infty}{\textrm{Law}}\mathcal N\left(0,\mu_2-4\mu_1^2(1-\alpha)^2\right).
	$
	\item[iii)] If $1/2<a<1$, then 
	$
	\sqrt{n^{2a-1}}\left(n^{1-a}\left(\frac{T_n}{n}-\frac{(1-\alpha)\mu _{1}}{1-a}\right)-L\right)\converge{12}{n}{+\infty}{\textrm{Law}}\mathcal N\left(0,\frac{\sigma^2}{2a-1}\right)
	$
			where $L$ is the random variable given in Theorem \ref{Loi_des_grands_nombres_sous_moment_ordre_deux}.
\end{itemize}
\end{Th}
\noindent\textbf{Remark 3}. Theorem \ref{TLC} is an extension of Theorems $3.3$ and $3.6$ in \cite{Bercu2018}, Theorem $1.1$ in \cite{Bertenghi2021} and Theorem $1.2$ in \cite{Bertoin2024}. Its proof is based on a martingale central limit theorem under the Lindeberg's condition (see section $2$ below).\\
\\
\noindent In the spirit of \cite{Gut--Stadtmuller2022} (see also \cite{Aguech--Elmachkouri2024} and \cite{Roy--Takei--Tanemura2024}), we investigate the asymptotic normality of the step-reinforced random walk in the gradually increasing memory setting. The idea is to introduce a gap (depending on $n$) in the memory of the random walk and to investigate its effect on the trajectories of the walker. More precisely, let $(m_n)_{n\geqslant 1}$ be a non-decreasing sequence of positive integers growing to infinity such that $m_n<n$ for any $n\geqslant 1$. Let $(X_n)_{n\geqslant 1}$ be the steps of the unbalanced step-reinforced random walk defined by \eqref{definition_unbalanced_step-reinforced_random_walk} and recall that $T_n=\sum_{k=1}^nX_k$ for any integer $n\geqslant 1$. For any $n\geqslant 1$ and any  $1\leqslant k\leqslant n$, we denote 
$$
\tilde{S}_{k,n}= \begin{cases}
	T_k & \textrm{if $\,1\leqslant k\leqslant m_n$}\\ 
	T_{m_n}+\sum_{j=m_n+1}^kX_j' & \textrm{if $\,m_n<k\leqslant n$}
\end{cases}
$$	
where
\begin{equation}\label{Def_X_prime}
	X_j'= \left\{
	\begin{array}{ccc}
		+X_{U_{j, n}} & \text { with probability } & p\alpha\\
		-X_{U_{j, n}} & \text { with probability } & (1-p)\alpha\\
		\xi_{j+1} & \text {with probability } & 1-\alpha
	\end{array}\right.
\end{equation}
for any $m_n<j \leq n$ and the random variables $\left\{U_{j, n}\,;m_n<j\leqslant n,\,n\geqslant 1\right\}$ are i.i.d and satisfy $U_{j,n}\hookrightarrow\mathcal U\left(\{1,...,m_n\}\right)$ for any $n\geqslant 1$ and any $m_n< j\leqslant n$. It is important to note that we no longer deal with time series but rather with triangular arrays of random variables introducing a gap in the memory of the random walk since any step $X'_j$ with $m_n<j\leqslant n$ does no longer depend on the steps $X_{m_n+1},\ldots,X_{j-1}$ but depends only on the steps $X_1,\ldots,X_{m_n}$.  In the sequel, we consider the unbalanced step-reinforced process with gradually increasing memory $(S_n)_{n\geqslant 1}$ defined by $S_0=0$ and $S_n=\tilde{S}_{n,n}$ for any integer $n\geqslant 1$ and we assume that $\lim_{n\to+\infty}n^{-1}m_n=\theta$ for some $\theta \in[0,1]$.
\begin{Th}\label{TLC_bis} Assume that $\mathbb E[\xi_1^2]<+\infty$ and denote $\tau=\theta+a(1-\theta)$.
\begin{itemize}
\item[i)] If $-1\leqslant a<1/2$, then $
		\sqrt{m_n}\left(\frac{S_n}{n}-\frac{(1-\alpha)\mu_1}{1-a}\right)\converge{12}{n}{+\infty}{\textrm{Law}}\mathcal N\left(0,\frac{\tau^2\sigma^2}{1-2a}+\theta(1-\theta)\sigma^2\right).
		$
\item[ii)] If $a=1/2$, then 
		$
		\frac{\sqrt{m_n}}{\sqrt{\log m_n}}\left(\frac{S_n}{n}-2(1-\alpha)\mu_1\right)\converge{12}{n}{+\infty}{\textrm{Law}}\mathcal N\left(0,\mu_2\tau^2-4\mu_1^2\tau^2(1-\alpha)^2\right).
		$
		\item[iii)] If $1/2<a<1$, then 
		$
		\sqrt{m_n^{2a-1}}\left(m_n^{1-a}\left(\frac{S_n}{n}-\frac{(1-\alpha)\mu_1}{1-a}\right)-\tau_nL\right)\converge{12}{n}{+\infty}{\textrm{Law}}\mathcal N\left(0,\frac{\tau^2\sigma^2}{2a-1}+\theta(1-\theta)\sigma^2\right)
		$
		where $\tau_n:=n^{-1}m_n+a\left(1-n^{-1}m_n\right)$ and $L$ is the random variable given in Theorem \ref{Loi_des_grands_nombres_sous_moment_ordre_deux}.
	\end{itemize}
\end{Th}
\noindent \textbf{Remark 4}. Theorem \ref{TLC_bis} is an extention of Theorem $3.1$ in \cite{Gut--Stadtmuller2022}, Theorem $1$ in \cite{Aguech--Elmachkouri2024_corrigendum} and Theorem $2.1$ in \cite{Roy--Takei--Tanemura2024}. Its proof could be carried on using Lindeberg's method following the approach in \cite{Aguech--Elmachkouri2024} but we opt for the characteristic function approach as in \cite{Roy--Takei--Tanemura2024}. The main challenge in the proof of Theorem \ref{TLC_bis} is to establish the almost sure convergence of some conditional characteristic functions with the constraint that the steps of the walker are no longer bounded.

\section{Proofs}
Let $(p_n)_{n\geqslant 1}$ and $(q_n)_{n\geqslant 1}$ be two arbitrary sequences of positive real numbers. In this section, the notation $p_n\lesssim q_n$ means that there exists a positive constant $C$ (not depending on $n$) such that $p_n\leqslant Cq_n$ for any integer $n\geqslant 1$ whereas the notation $p_n\sim_{n\to+\infty}q_n$ means that $\lim_{n\to+\infty}p_n/q_n=1$. As usual, for any $q\geqslant 1$, the notation $\Vert\,.\,\Vert_q$ stands for the standard $\mathbb L^q$-norm with respect to the probability measure $\P$.   Now, one can notice that the unbalanced step-reinforced $(X_n)_{n\geqslant 1}$ with parameters $p\in [0,1]$ and $\alpha\in[0,1]$ given by \eqref{definition_unbalanced_step-reinforced_random_walk} can be defined for any integer $n\geqslant 1$ by
\begin{equation}\label{def_sympa_step-reinforced_random_walk}
X_{n+1}=\gamma_n\rho_nX_{U_n}+(1-\gamma_n)\xi_{n+1}
\end{equation}
where $\gamma_n$ and $\rho_n$ are two independent random variables defined by
$ \P(\gamma_n=1)=\alpha=1-\P(\gamma_n=0)$ and $\P(\rho_n=1)=p=1-\P(\rho_n=-1)$ such that $\gamma_n$ and $\rho_n$ are independent of $U_1,...,U_n, X_1,...,X_n, \xi_1,\xi_2,...$. If $Z$ is any integrable random variable and $\mathcal G$ is a $\sigma$-algebra, we denote by $\mathbb E_{\mathcal G}[Z]$ and $\mathbb V_{\mathcal G}[Z]$ the conditional expectation and the conditional variance of $Z$ with respect to the $\sigma$-algebra $\mathcal G$ respectively. Using \eqref{def_sympa_step-reinforced_random_walk} and considering the filtration $(\F_k)_{k\geqslant 0}$ defined by $\mathcal F_0=\sigma(\emptyset, \Omega)$ and $\mathcal F_k=\sigma\left(\xi_j,\,\rho_j,\,U_j,\,\gamma_j;\,j\leqslant k\right)$ for any integer $k\geqslant 1$, we have
$$
\mathbb E_{\F_n}[X_{n+1}]=\frac{a}{n}\sum_{k=1}^nX_k+(1-\alpha)\mu_1.
$$
Consequently,
$$
\mathbb E_{\F_n}\!\left[\sum_{k=1}^{n+1}X_k\right]=\left(1+\frac{a}{n}\right)\sum_{k=1}^nX_k+(1-\alpha)\mu_1.
$$
So, 
$$
\mathbb E_{\F_n}[\tilde{T}_{n+1}]=\left(1+\frac{a}{n}\right)\tilde{T}_n\quad\textrm{where}\quad \tilde{T}_n:=\sum_{k=1}^n\left(X_k-\frac{(1-\alpha)\mu_1}{1-a}\right).
$$
Following the martingale approach in \cite{Bercu2018}, we define the martingale $(M_n)_{n\geqslant 1}$ with respect to the filtration $(\F_n)_{n\geqslant 1}$ by 
\begin{equation}\label{def_M_n}
M_n=a_n\tilde{T}_n=a_n\sum_{k=1}^n\left(X_k-\frac{(1-\alpha)\mu_1}{1-a}\right)
\end{equation}
with $a_1=1$ and for any integer $n\geqslant 2$,
\begin{equation}\label{def_a_n}
	a_n=\left\{\begin{array}{ccc}
		\frac{\Gamma(n)\Gamma(a+1)}{\Gamma(n+a)} & \textrm{if} & -1<a<1\\[1mm]
		2(n-1) & \textrm{if} & a=-1
	\end{array}\right.
\end{equation}
where $\Gamma$ stands for the Euler Gamma function. If we denote $A_n:=\sum_{k=1}^na_k$ and $v_n:=\sum_{k=1}^na_k^2$, then for $a=-1$, we have $
A_n=n^2-n+1$ and $v_n=1+\frac{2}{3}(n-1)n(2n-1)$. So, with the convention $\Gamma(0):=2$, for any $-1\leqslant a<1$, we have 
\begin{equation}\label{estimee_de_a_n}
	\lim_{n\to+\infty}\frac{n^{a}a_n}{\Gamma(a+1)}=1,\qquad\lim_{n\to+\infty}\frac{(1-a)A_n}{n^{1-a}\Gamma(a+1)}=1
\end{equation}
and
\begin{equation}\label{estimee_de_v_n}
	\lim_{n\to+\infty}r_n^{-1}v_n=1\quad\textrm{with}\quad r_n:=\left\{\begin{array}{ccc}
		\frac{n^{1-2a}\Gamma^2(a+1)}{1-2a} & \textrm{if} & -1\leqslant a<1/2\\[1mm]
		\frac{\pi}{4}\log n & \textrm{if} &a=1/2\\[1mm]
		\sum_{k=1}^{\infty}a^2_k<+\infty &\textrm{if}& 1/2<a<1.
	\end{array}\right.
\end{equation}
In the sequel, the following lemma will be useful.
\begin{Lem}\label{Lemme_technique}
	If $\psi :\R\to\R$ is a measurable function such that $\psi(-x)=\psi(x)$ for any $x\in\R$, then $\mathbb E[\psi(X_k)]=\mathbb E[\psi\left(\xi_1\right)]$ for any integer $k\geqslant 1$.
\end{Lem}
\noindent {\em Proof of Lemma \ref{Lemme_technique}}. Since $X_1=\xi_1$ a.s. the result is true for $k=1$. Now, let $k\geqslant 1$ and assume that $\mathbb E[\psi(X_j)]=\mathbb E[\psi\left(\xi_j\right)]$ for any $1\leqslant j\leqslant k$. Since $X_{k+1}=\gamma_k\rho_kX_{U_k}+(1-\gamma_k)\xi_{k+1}$, we derive
\begin{align*}
	\mathbb E[\psi\left(X_{k+1}\right)]
	&=\mathbb E[\psi(X_{k+1})\vert\gamma_k=0]\P(\gamma_k=0)+\mathbb E[\psi(X_{k+1})\vert\gamma_k=1]\P(\gamma_k=1)\\
	&=(1-\alpha)\mathbb E[\psi(\xi_{k+1})]+\alpha\mathbb E[\psi(X_{U_k})]\quad\textrm{since}\quad \psi(\rho_kX_{U_k})=\psi(X_{U_k})\,\, \textrm{a.s.}\\[-1mm]
	&=(1-\alpha)\mathbb E[\psi(\xi_{k+1})]+\alpha\sum_{\ell=1}^k\mathbb E[\psi(\xi_{\ell})]\P(U_k=\ell)\\
	&=(1-\alpha)\mathbb E[\psi(\xi_{k+1})]+\alpha\mathbb E[\psi(\xi_{k+1})]=\mathbb E[\psi(\xi_1)].
\end{align*}
The proof of Lemma \ref{Lemme_technique} is complete.$\hfill\Box$\\
\\
Using \eqref{def_M_n}, we define 
$$
M_0:=\mathbb E[M_1]=\mu_1-\frac{(1-\alpha)\mu_1}{1-a}=\frac{(\alpha-a)\mu_1}{1-a}.
$$
Consequently, for any $n\geqslant 1$, 
$$
M_n=M_0+\sum_{k=1}^n\left(M_k-M_{k-1}\right)=\frac{(\alpha-a)\mu_1}{1-a}+\sum_{k=1}^na_k\varepsilon_k
$$
with $\varepsilon_k=X_k-\mathbb E_{\F_{k-1}}\left[X_k\right]=\tilde{\varepsilon}_k+\overline{\varepsilon}_k$ where
\begin{equation}\label{def_epsilon_barre_k}
\tilde{\varepsilon}_k:=X_k\ind_{\vert X_k\vert\leqslant c_k}-\mathbb E_{\F_{k-1}}[X_k\ind_{\vert X_k\vert\leqslant c_k}]\quad\textrm{and}\quad
\overline{\varepsilon}_k
:=X_k\ind_{\vert X_k\vert> c_k}-\mathbb E_{\F_{k-1}}[X_k\ind_{\vert X_k\vert>c_k}]
\end{equation}
and $(c_k)_{k\geqslant 1}$ is an increasing sequence of positive real numbers satisfying $\lim_{k\to+\infty}c_k=+\infty$.
\begin{Lem}\label{Convergence_de_moment_ordre_2_de_epsilon_barre_k}
	If $-1\leqslant a<1$ and $\mathbb E[\xi_1^2]<+\infty$, then
	\begin{itemize}
		\item[i)] $n^{-1}\sum_{j=1}^nX_j\converge{19}{n}{+\infty}{$\mathbb L^2$ and a.s.}\frac{(1-\alpha)\mu_1}{1-a}$ and $n^{-1}\sum_{j=1}^nX_j^2\converge{19}{n}{+\infty}{$\mathbb L^1$ and a.s.}\mu_2$.
		\item[ii)] $\mathbb E_{\F_{k-1}}[X_k]\converge{17}{k}{+\infty}{$\mathbb L^2$ and a.s.}\frac{(1-\alpha)\mu_1}{1-a}$,\,\, $\mathbb E_{\F_{k-1}}[X_k^2]\converge{17}{k}{+\infty}{$\mathbb L^1$ and a.s.}\mu_2$ and $\mathbb E_{\F_{k-1}}[\varepsilon_k^2]\converge{17}{k}{+\infty}{$\mathbb L^1$ and a.s.}\sigma^2$ where $\sigma^2$ is defined by \eqref{definition_sigma_carre_et_a}.
		\item[iii)] $\mathbb E_{\F_{k-1}}[\tilde{\varepsilon}_k^2]\converge{8}{k}{+\infty}{$\mathbb L^1$}\sigma^2$ and $\lim_{k\to+\infty}\mathbb E\left[\overline{\varepsilon}^2_k\right]=0$.
	\end{itemize}
\end{Lem}
\noindent {\em Proof of Lemma \ref{Convergence_de_moment_ordre_2_de_epsilon_barre_k}}. 
Let $n\geqslant 1$ be fixed and recall that $M_n=a_n\tilde{T}_n=\frac{(\alpha-a) \mu_1}{1-a}+\sum_{k=1}^na_k\varepsilon_k$ where $\tilde{T}_n=\sum_{k=1}^n\left(X_k-\frac{(1-\alpha)\mu_1}{1-a}\right)$ and $\varepsilon_k=\tilde{\varepsilon}_k+\overline{\varepsilon}_k=X_k-\mathbb E_{\F_{k-1}}[X_k]$. Using Lemma \ref{Lemme_technique}, we have $\mathbb E[\varepsilon_k^2]\leqslant \mathbb E[\xi_1^2]$ for any integer $k\geqslant 1$. Moreover, using the inequality $(a+b)^2\leqslant 2(a^2+b^2)$ for any $(a,b)\in\R^2$, we derive
$$
\frac{\mathbb E[\tilde{T}_n^2]}{n^2}=\frac{\mathbb E[M_n^2]}{(na_n)^2}\leqslant \frac{2}{(na_n)^2}\left(\frac{(\alpha-a)^2\mu_1^2}{(1-a)^2}+\sum_{k=1}^na_k^2\mathbb E[\varepsilon_k^2]\right)\leqslant \frac{2}{(na_n)^2}\left(\frac{(\alpha-a)^2\mu_1^2}{(1-a)^2}+v_n\mathbb E[\xi_1^2]\right).
$$
Using \eqref{def_a_n}, \eqref{estimee_de_a_n} and  \eqref{estimee_de_v_n}, we obtain $n^{-1}\tilde{T}_n\converge{10}{n}{+\infty}{$\mathbb L^2$}0$, that is 
\begin{equation}\label{loi_des_grands_nombres_dans_L_2}
	\frac{1}{n}\sum_{k=1}^nX_k\converge{10}{n}{+\infty}{$\mathbb L^2$}\frac{(1-\alpha)\mu_1}{1-a}.
\end{equation}
Since 
$
\mathbb E_{\F_n}[\tilde{T}_{n+1}]=\left(1+\frac{a}{n}\right)\tilde{T}_n,
$ 
we have
$
\mathbb E_{\F_n}\left[\frac{\tilde{T}_{n+1}}{n+1}\right]-\frac{\tilde{T}_n}{n}=\left(\frac{a-1}{n+1}\right)\frac{\tilde{T}_n}{n}.
$
Using \eqref{estimee_de_a_n}, we obtain 
$$
\left\Vert\mathbb E_{\F_n}\!\left[\frac{\tilde{T}_{n+1}}{n+1}\right]-\frac{\tilde{T}_n}{n}
\right\Vert_1\leqslant\frac{2\Vert M_n\Vert_2}{n(n+1)a_n}\lesssim\frac{\frac{(\alpha-a)\mu_1}{1-a}+\sqrt{\sum_{j=1}^na_j^2\mathbb E[\varepsilon_j^2]}}{n^{2-a}}\leqslant\frac{\frac{(\alpha-a)\mu_1}{1-a}+\sqrt{v_n\mu_2}}{n^{2-a}}.
$$
Using \eqref{estimee_de_v_n}, we derive that $\left(k^{-1}\tilde{T}_k\right)_{k\geqslant 1}$ is a quasi-martingale in the sense that 
$$
\sum_{n\geqslant 1}\left\Vert\mathbb E_{\F_n}\!\left[\frac{\tilde{T}_{n+1}}{n+1}\right]-\frac{\tilde{T}_n}{n}
\right\Vert_1<+\infty\quad\textrm{and}\quad 
\sup_{k\geqslant 1}\frac{\mathbb E[\vert \tilde{T}_k\vert]}{k}\leqslant \mathbb E\vert \xi_1\vert+\frac{(1-\alpha)\mu_1}{1-a}<+\infty.
$$
Using Theorem 9.4 in \cite{Metivier1982} (see also Theorem 1.16 in \cite{Qin2024}), there exists an integrable random variable $Z_{\infty}$ such that $k^{-1}\tilde{T}_k\converge{10}{k}{+\infty}{a.s.}Z_{\infty}$ and 
$
\mathbb E\vert Z_{\infty}\vert\leqslant \liminf_{k\to+\infty}\mathbb E\left\vert k^{-1}\tilde{T}_k\right\vert=0,
$ that is $Z_{\infty}=0$ a.s. which means that 
\begin{equation}\label{loi_forte_des_grands_nombres}
	\frac{1}{n}\sum_{k=1}^nX_k\converge{12}{n}{+\infty}{a.s.}\frac{(1-\alpha)\mu_1}{1-a}.
\end{equation}
Combining \eqref{loi_des_grands_nombres_dans_L_2} and \eqref{loi_forte_des_grands_nombres}, we get
\begin{equation}\label{convergence_esperance_conditionnelle_X_k_dans_L_deux}
	\mathbb E_{\F_{k-1}}[X_k]=\frac{a}{k-1}\sum_{\ell=1}^{k-1}X_{\ell}+(1-\alpha)\mu_1\converge{19}{k}{+\infty}{$\mathbb L^2$ and a.s.}\frac{(1-\alpha)\mu_1}{1-a}.
\end{equation}
Using \eqref{def_sympa_step-reinforced_random_walk}, we have $X_1^2=\xi_1^2$ a.s. and $X_{k+1}^2=\gamma_kX_{U_k}^2+(1-\gamma_k)\xi_{k+1}^2$ a.s. for any integer $k\geqslant 1$. It means that $(X_k^2)_{k\geqslant 1}$ is a positively step-reinforced random walk (see equation $(1.1)$ in \cite{Hu--Zhang2024}). So, using Theorem 1.1 in \cite{Hu--Zhang2024}, we have $n^{-1}\sum_{k=1}^n\left(X_k^2-\mu_2\right)\converge{12}{n}{+\infty}{a.s.}0$. Using Lemma \ref{Lemme_technique}, we have $\lim_{c\to+\infty}\sup_{k\geqslant 1}\mathbb E[X^2_k\ind_{\vert X_k\vert>c}]=\lim_{c\to+\infty}\mathbb E[\xi^2_1\ind_{\vert \xi_1\vert>c}]=0$ which means that the random variables $(X_k^2)_{k\geqslant 1}$ are uniformly integrable. So, the random variables $(n^{-1}\sum_{j=1}^n(X_j^2-\mu_2))_{n\geqslant 1}$ are also uniformly integrable and consequently
%\begin{equation}\label{LLN_forte_for_X_k_carree}
$
n^{-1}\sum_{k=1}^n\left(X_k^2-\mu_2\right)\converge{19}{n}{+\infty}{$\mathbb L^1$ and a.s.}0.	
$	
%\end{equation}
Using the fact that $\mathbb E_{\F_{k-1}}[X_k^2]-\mu_2=\frac{\alpha}{k-1}\sum_{j=1}^{k-1}\left(X_j^2-\mu_2\right)$ for any integer $k\geqslant 2$, we get
\begin{equation}\label{convergence_L_deux_esp_cond_X_k_carree}
	\mathbb E_{\F_{k-1}}[X_k^2]\converge{19}{k}{+\infty}{\textrm{$\mathbb L^1$ and a.s.}}\mu_2.
\end{equation} 
Combining \eqref{convergence_esperance_conditionnelle_X_k_dans_L_deux},  \eqref{convergence_L_deux_esp_cond_X_k_carree} and noting that 
$$
\left\Vert \left(\mathbb E_{\F_{k-1}}[X_k]\right)^2-\tau^2\right\Vert_1\leqslant \left\Vert \mathbb E_{\F_{k-1}}[X_k]-\tau\right\Vert_2^2+2\tau\left\Vert\mathbb E_{\F_{k-1}}[X_k]-\tau\right\Vert_1\converge{12}{k}{+\infty}{}0
$$
where $\tau:=\frac{(1-\alpha)\mu_1}{1-a}$, we obtain 
$$
\mathbb E_{\F_{k-1}}[\varepsilon_k^2]=\mathbb E_{\F_{k-1}}[X_k^2]-\left(\mathbb E_{\F_{k-1}}[X_k]\right)^2\converge{19}{k}{+\infty}{\textrm{$\mathbb L^1$ and a.s.}}\mu_2-\frac{(1-\alpha)^2\mu_1^2}{(1-a)^2}=\sigma^2.
$$
Using Lemma \ref{Lemme_technique}, we have $\mathbb E[\overline{\varepsilon}_k^2]\leqslant \mathbb E[\xi_1^2\ind_{\vert \xi_1\vert>c_k}]$ and we get $\lim_{k\to+\infty}\mathbb E[\overline{\varepsilon}_k^2]=0$.
Finally, for any integer $k\geqslant 1$,
%\begin{equation}\label{decomposition_esp_cond_epsilon_tild_carre}
$$	
\mathbb E_{\F_{k-1}}[\tilde{\varepsilon}_k^2]
	=\mathbb E_{\F_{k-1}}[\varepsilon_k^2]-\mathbb E_{\F_{k-1}}[\overline{\varepsilon}_k^2]
	-2\mathbb E_{\F_{k-1}}[\tilde{\varepsilon}_k\overline{\varepsilon}_k].
%\end{equation}
$$
Noting that 
$
\mathbb E[\vert\mathbb E_{\F_{k-1}}[\tilde{\varepsilon}_k\overline{\varepsilon}_k]\vert]\leqslant\sqrt{\mathbb E[\tilde{\varepsilon}_k^2]\mathbb E[\overline{\varepsilon}_k^2]}\leqslant \sqrt{\mu_2\mathbb E[\overline{\varepsilon}_k^2]}\converge{10}{k}{+\infty}{}0
$, we get $\mathbb E_{\F_{k-1}}[\tilde{\varepsilon}_k^2]\converge{10}{k}{+\infty}{\textrm{$\mathbb L^1$}}\sigma^2$. The Proof of Lemma \ref{Convergence_de_moment_ordre_2_de_epsilon_barre_k} is complete.$\hfill\Box$
\begin{Lem}\label{Lemme_de_l_algorithme}
	Let $(s_k)_{k\geqslant 1}$ be a sequence of real numbers such that $\lim_{k\to+\infty}s_k=s\in\R$. Let also $x_1\in\R$ and $-1\leqslant t<1$ be fixed and define $x_{k+1}=s_k+\frac{t}{k}\sum_{j=1}^kx_j$ for any integer $k\geqslant 1$. Then, the sequence $(x_k)_{k\geqslant 1}$ is convergent and $\lim_{k\to+\infty}x_k=\lim_{k\to+\infty}\frac{1}{k}\sum_{j=1}^kx_j=s/(1-t)$.
\end{Lem}
\noindent {\em Proof of Lemma \ref{Lemme_de_l_algorithme}}. 
Note that $V_{k+1}=\left(1+\frac{t}{k}\right)V_k+s_k$ where $V_k=\sum_{j=1}^kx_j$ for any $k\geqslant 1$. So, for any $k\geqslant 2$, we get 
$$
\frac{V_k}{k}=\frac{x_1\Gamma(k+t)}{k\Gamma(k)\Gamma(t+1)}+\frac{\Gamma(k+t)}{\Gamma(k+1)}\sum_{j=1}^{k-1}\frac{\Gamma(j+1)s_j}{\Gamma(j+t+1)}\sim_{k\to+\infty}\frac{x_1k^{t-1}}{\Gamma(t+1)}+\frac{1}{k^{1-t}}\sum_{j=1}^{k-1}\frac{s_j}{j^t}\converge{12}{k}{+\infty}{}\frac{s}{1-t}.
$$
The proof of Lemma \ref{Lemme_de_l_algorithme} is complete.$\hfill\Box$\\
\\
{\em Proof of Theorem \ref{Loi_forte_des_grands_nombres}}. 
Using Lemma \ref{Lemme_technique}, we have $$
\sum_{k=1}^{+\infty} \mathbb{P}\left(\left|X_k\right|>k\right)
=\sum_{k=1}^{+\infty} \mathbb{P}\left(\vert\xi_1\vert>k\right)\leqslant \mathbb E[\vert \xi_1\vert]<+\infty.
$$ 
Using Borel Cantelli's lemma, we obtain $\mathbb{P}( \overline{N}) = 1$ where $\overline{N}:=\cup_{n\geqslant 1} \cap_{k \geqslant n}\{\vert X_k\vert \leqslant k\}$. 
So, if $\omega\in\overline{N} $, there exists $K=K(\omega) \in \mathbb{N} $ such that $\vert X_k(\omega)\vert \leqslant k$ for any integer $k\geqslant K$. Consequently, for any $n\geqslant K$, we have
\begin{equation}\label{convergence_du_reste}
\frac{1}{n} \sum_{k=1}^n \vert X_k(\omega)\vert \ind_{|X_k(\omega)|>k}  
= \frac{1}{n} \sum_{k=1}^{K} \vert X_k(\omega)\vert \ind_{|X_k(\omega)|>k}  
\xrightarrow[n \to +\infty]{} 0.
\end{equation}
So, it is sufficient to prove that 
\begin{equation}\label{convergence_moyennes_phi_k_de_X_k}
\frac{1}{n} \sum_{k=1}^n \varphi_k(X_k) \converge{12}{n}{+\infty}{a.s.}
\frac{(1-\alpha) \mu_1}{1-a}
\end{equation} 
where $\varphi_k$ is the function defined by $\varphi_k(x)=x\ind_{\vert x\vert \leq k}$ for any $x\in\R$ and any $k\geqslant 1$. In order to obtain \eqref{convergence_moyennes_phi_k_de_X_k}, it is sufficient to prove that 
\begin{equation}\label{convergence_moyennes_esperance_conditionnelle_de_phi_k_de_X_k}
\frac{1}{n} \sum_{k=1}^n \mathbb{E}_{\F_{k-1}}[\varphi_k(X_k)]\converge{12}{n}{+\infty}{a.s.}\frac{(1-\alpha) \mu_1}{1-a}.
\end{equation}
In fact, if $Y_k:=\varphi_k(X_k)-\mathbb{E}_{\F_{k-1}}[\varphi_k(X_k)]$, then 
$$
\sum_{k=1}^{+\infty} \frac{\mathbb{E}[Y_k^2]}{k^2}\leqslant \sum_{k=1}^{+\infty} \frac{\mathbb{E}[X_k^2\ind_{\vert X_k\vert \leqslant k}]}{k^2} = \sum_{k=1}^{+\infty} \frac{\mathbb{E}[\xi_1^2\ind_{\vert \xi_1\vert \leqslant k}]}{k^2}\leqslant \mathbb{E}\!\!\left[\xi_1^2 \sum_{k\geqslant \vert\xi_1\vert\vee 1}\frac{1}{k^2}\right] \leqslant \mathbb{E}[\vert\xi_1\vert]<+\infty.
$$
Consequently, we get $\sum_{k\geqslant 1}k^{-2}\mathbb E_{\F_{k-1}}[Y_k^2]<+\infty$ a.s. Noting that $(Y_k/k)_{k\geqslant 1}$ is a sequence of martingale differences and applying Theorem 2.15 in \cite{Hall--Heyde1980}, we derive that $\sum_{k\geqslant 1}k^{-1}Y_k$ converges a.s. Using Kronecker's lemma, we get 
\begin{equation}\label{convergence_moyennes_des_Y_k}
\frac{1}{n}\sum_{k=1}^n Y_k\converge{12}{n}{+\infty}{a.s.}0.
\end{equation}
Now, in order to prove \eqref{convergence_moyennes_esperance_conditionnelle_de_phi_k_de_X_k}, we use \eqref{def_sympa_step-reinforced_random_walk} and we get
$$
\mathbb{E}_{\F_k}[\varphi_{k+1}(X_{k+1})]=(1-\alpha) \mathbb{E}[\varphi_{k+1}(\xi_1)] + \frac{a}{k} \displaystyle\sum_{j=1}^k \varphi_{k+1}\left(X_j\right).
$$
Denoting 
$$
Z_k:=\frac{a}{k}\sum_{j=1}^k\left(\varphi_{k+1}(X_j)-\varphi_j(X_j)\right)+\frac{a}{k}\sum_{j=1}^kY_j,
$$
we get 
\begin{equation}\label{equation_clef}
\mathbb{E}_{\F_k}[\varphi_{k+1}(X_{k+1})]=Z_k+(1-\alpha) \mathbb{E}[\varphi_{k+1}(\xi_1)] + \frac{a}{k} \displaystyle\sum_{j=1}^k\mathbb{E}_{\F_{j-1}}[\varphi_{j}(X_{j})].
\end{equation}
Moreover, using \eqref{convergence_du_reste} and \eqref{convergence_moyennes_des_Y_k}, we have
\begin{equation}\label{Majoration_Z_k}
\vert Z_k\vert\leqslant \frac{\vert a\vert}{k}\sum_{j=1}^k\vert X_j\vert\ind_{\vert X_j\vert>j}+\left\vert\frac{a}{k} \sum_{j=1}^kY_j\right\vert\converge{12}{k}{+\infty}{a.s.}0.
\end{equation}
Combining \eqref{equation_clef}, \eqref{Majoration_Z_k} and Lemma \ref{Lemme_de_l_algorithme}, we obtain 
\begin{equation}\label{convergence_moyennes_conditionnelles_des_phi_k_X_k}
\frac{1}{k}\sum_{j=1}^k\mathbb{E}_{\F_{j-1}}[\varphi_{j}(X_{j})]\converge{12}{k}{+\infty}{a.s.}\frac{(1-\alpha)\mu_1}{1-a}.
\end{equation}
Combining \eqref{convergence_moyennes_des_Y_k} and \eqref{convergence_moyennes_conditionnelles_des_phi_k_X_k}, we derive \eqref{convergence_moyennes_phi_k_de_X_k}. The proof of Theorem \ref{Loi_forte_des_grands_nombres} is complete.$\hfill\Box$\\
\\
{\em Proof of Theorem \ref{Loi_des_grands_nombres_Marcinkiewicz-Zygmund}}. Let $n\geqslant 1$ be a fixed integer and recall the notations:
$$
\tilde{T}_n:=\sum_{k=1}^n\left(X_k-\frac{(1-\alpha)\mu_1}{1-a}\right)\quad\textrm{and}\quad M_n=a_n\tilde{T}_n=\frac{(\alpha-a)\mu_1}{1-a}+\sum_{k=1}^na_{k}\varepsilon_{k}
$$ 
where $\varepsilon_{k}=X_{k}-\mathbb E_{\F_{k-1}}[X_{k}].$
Let $M_n^{\ast}:=\max_{1\leqslant k\leqslant n}\vert M_k-\frac{(\alpha-a) \mu_1}{1-a}\vert$ and $1<r<2$ and denote $d_n:=n^{\frac{1}{r}-\frac{1}{2}}\sqrt{v_n}$ where $v_n=\sum_{k=1}^na_k^2$. Our aim is to prove that $\sum_{n=1}^{\infty}  \frac{1}{n}\mathbb{P}\left( M_n^{\ast}>\eta d_n\right)<+\infty$ for any $\eta>0$. Let $(c'_n)_{n\geqslant 1}$ be the increasing sequence of positive real numbers defined by $c'_n=n^{1/r}$ for any integer $n\geqslant 1$ and denote
$$
\tilde{\varepsilon}_k:=X_k\ind_{\vert X_k\vert\leqslant c'_n}-\mathbb E_{\F_{k-1}}[X_k\ind_{\vert X_k\vert\leqslant c'_n}]\quad\textrm{and}\quad
\overline{\varepsilon}_k
:=X_k\ind_{\vert X_k\vert> c'_n}-\mathbb E_{\F_{k-1}}[X_k\ind_{\vert X_k\vert>c'_n}]
$$
for any $1\leqslant k\leqslant n$, so that  $M_n=\frac{(\alpha-a)\mu_1}{1-a}+\tilde{M}_n+\overline{M}_n$ with $\tilde{M}_n:=\sum_{k=1}^na_k\tilde{\varepsilon}_k$ and $\overline{M}_n:=\sum_{k=1}^na_k\overline{\varepsilon}_k$. Let $\eta>0$ be fixed. Then, 
$$
\P\left(M_n^{\ast}>2\eta d_n\right)\leqslant
\P\left(\tilde{M}_n^{\ast}>\eta d_n\right)+\P\left(\overline{M}_n^{\ast}>\eta d_n\right)
$$
with $\tilde{M}_n^{\ast}:=\max_{1\leqslant k\leqslant n}\vert \tilde{M}_k\vert$ and $\overline{M}_n^{\ast}:=\max_{1\leqslant k\leqslant n}\vert \overline{M}_k\vert$. Moreover, using the Doob maximal inequality for martingales (see \cite[Theorem 2.11]{Hall--Heyde1980}), there exists a positive constant $\kappa>0$ such that
$$
\mathbb E[\tilde{M}^{\ast2}_n]\leqslant \kappa \sum_{k=1}^na_k^2\mathbb E[\tilde{\varepsilon}_k^2].
$$
Using Markov's inequality, we get
$$
\P\left(\tilde{M}_n^{\ast}>\eta d_n\right)\leqslant \frac{\kappa\sum_{k=1}^na_k^2\mathbb E[\tilde{\varepsilon}_k^2]}{\eta^2d_n^2}\leqslant \frac{\kappa\sum_{k=1}^na_k^2\mathbb E[X_k^2\ind_{\vert X_k\vert\leqslant c'_n}]}{\eta^2d_n^2}.
$$
Using Lemma \ref{Lemme_technique}, we have $\mathbb E[X_k^2\ind_{\vert X_k\vert\leqslant c'_n}]=\mathbb E[\xi_1^2\ind_{\vert \xi_1\vert\leqslant c'_n}]$. So, for $-1\leqslant a<1$, we derive 
\begin{align*}
\sum_{n=1}^{+\infty}\frac{1}{n}\P\left(\tilde{M}_n^{\ast}>\eta d_n\right)&\leqslant\frac{\kappa}{\eta^2} \sum_{n=1}^{+\infty}\frac{v_n}{nd_n^2}\mathbb E[\xi_1^2\ind_{\vert \xi_1\vert\leqslant c'_n}]=\frac{\kappa}{\eta^2} \sum_{n=1}^{+\infty}n^{-\frac{2}{r}}\mathbb E[\xi_1^2\ind_{\vert \xi_1\vert\leqslant n^{1/r}}]\\
&\leqslant \frac{\kappa}{\eta^2} \mathbb E\left[\xi_1^2\sum_{n\geqslant \vert\xi_1\vert^r}n^{-\frac{2}{r}}\right]\lesssim \frac{\mathbb E[\vert\xi_1\vert^r]}{\frac{2}{r}-1}<+\infty.
\end{align*}
Moreover, keeping in mind that $A_n:=\sum_{k=1}^na_k$, we have 
$$
\P\left(\overline{M}_n^{\ast}>\eta d_n\right)\leqslant \frac{\mathbb E[\vert\overline{M}_n^{\ast}\vert]}{\eta d_n}\leqslant \frac{1}{\eta d_n}\sum_{k=1}^na_k\mathbb E[\vert\overline{\varepsilon}_k\vert]\leqslant\frac{2A_n\mathbb E[\vert \xi_1\vert\ind_{\vert \xi_1\vert>c'_n}]}{\eta n^{\frac{1}{r}-\frac{1}{2}}\sqrt{v_n}}.
$$
Noting that $A_n\leqslant\sqrt{nv_n}$, we derive 
$$
\sum_{n=1}^{+\infty}\frac{1}{n}\P\left(\overline{M}_n^{\ast}>\eta d_n\right)\lesssim \frac{2}{\eta}\sum_{n=1}^{+\infty}n^{-1/r}\mathbb E[\vert \xi_1\vert\ind_{\vert\xi_1\vert>n^{1/r}}].
$$
Since $\mathbb E[\vert \xi_1\vert\ind_{\vert\xi_1\vert>n^{1/r}}]=\sum_{k=n}^{+\infty}\mathbb E\left[\vert\xi_1\vert\ind_{k<\vert\xi_1\vert^r\leqslant k+1}\right]$, we get 
\begin{align*}
\sum_{n=1}^{+\infty}\frac{1}{n}\P\left(\overline{M}_n^{\ast}>\eta d_n\right)&\lesssim \frac{2}{\eta}\sum_{k=1}^{+\infty}\left(\sum_{n=1}^{k}n^{-1/r}\right)\mathbb E[\vert\xi_1\vert\ind_{k<\vert\xi_1\vert^r\leqslant k+1}]\lesssim \frac{1}{1-\frac{1}{r}}\sum_{k=1}^{+\infty}k^{1-\frac{1}{r}}\mathbb E[\left\vert\xi_1\right\vert\ind_{k<\vert\xi_1\vert^r\leqslant k+1}]\\
&\leqslant \frac{1}{1-\frac{1}{r}}\sum_{k=1}^{+\infty}\mathbb E[\left\vert\xi_1\right\vert^r\ind_{k<\vert\xi_1\vert^r\leqslant k+1}]\leqslant\frac{\mathbb E[\vert\xi_1\vert^r]}{1-\frac{1}{r}}<+\infty.
\end{align*}
Finally, for any $\eta>0$ and any $-1\leqslant a<1$, we proved that
\begin{equation}\label{serie_convergente_fondamentale}
\sum_{n=1}^{+\infty}\frac{1}{n}\P\left(M_n^{\ast}>\eta d_n\right)<+\infty.
\end{equation} 
Now, one can notice that $(M_n^{*})_{n\geqslant 1}$ and $(d_n)_{n\geqslant 1}$ are non-decreasing sequences. So, for any integer $k\geqslant 1$ and any $2^k \leqslant n < 2^{k+1}$, we have
$$
\P(M_{2^k}^*> \eta d_{2^{k+1}})
\leqslant
\P(M_n^*>\eta d_n).
$$
Thus
$$
\frac{2^k}{2^{k+1}}
\P(M_{2^k}^* > \eta d_{2^{k+1}})
\leqslant
\sum_{2^k \leqslant n < 2^{k+1}}
\frac{1}{n}\P(M_n^* > \eta d_n).
$$
Hence
$$
\frac{1}{2}
\sum_{k=1}^{+\infty}
\P(M_{2^k}^* > \eta d_{2^{k+1}})
\leqslant
\sum_{k=1}^{+\infty}
\sum_{2^k\leqslant n < 2^{k+1}}
\frac{1}{n}\P(M_n^* > \eta\,d_n)\leqslant\sum_{n=1}^{+\infty}\frac{1}{n}\P\left(M_n^{\ast}>\eta d_n\right).
$$
Using \eqref{serie_convergente_fondamentale}, we obtain $\sum_{k=1}^{+\infty}
\P(M_{2^k}^* > \eta d_{2^{k+1}})<+\infty$ for any $\eta>0$. So, 
using Borel-Cantelli's lemma, we get  $d_{2^{k+1}}^{-1}M_{2^k}^*\converge{10}{k}{+\infty}{a.s.}0$. Consequently, since 
$
d_{2^{k+1}}^{-1}M_{2^k}^*
\leqslant 
d_n^{-1}M_n^*
\leqslant 
d_{2^k}^{-1}M_{2^{k+1}}^*
$ a.s. for any integer $k\geqslant 1$ and any $2^k \le n < 2^{k+1}$, we get $d_n^{-1}M_n^{*}\converge{12}{n}{+\infty}{a.s.}0$. In particular, we have $d_n^{-1}M_n\converge{12}{n}{+\infty}{a.s.}0$ which means that 
$$
\dfrac{na_n}{n^{\frac{1}{r}-\frac{1}{2}}\sqrt{v_n}} \left(\frac{T_n}{n}-\frac{(1-\alpha)\mu_1}{1-a}\right) \xrightarrow[n \to +\infty]{a.s}0.
$$
Using \eqref{estimee_de_a_n} and \eqref{estimee_de_v_n}, we get the desired results. The proof of Theorem \ref{Loi_des_grands_nombres_Marcinkiewicz-Zygmund} is complete.$\hfill\Box$\\
\\
{\em Proof of Theorem  \ref{Loi_des_grands_nombres_sous_moment_ordre_deux}}. Let $n\geqslant 1$ be a fixed integer. Following the martingale approach in \cite{Bercu2018}, we consider the predictable quadratic variation $\langle M\rangle_n$ of $M_n$ given by $\langle M\rangle_0=0$ and $\langle M\rangle_n=\sum_{k=1}^na_k^2\mathbb E_{\F_{k-1}}[\varepsilon_k^2]$. Using Lemma \ref{Convergence_de_moment_ordre_2_de_epsilon_barre_k}, we have $\langle M\rangle_n\sim_{n\to+\infty}\sigma^2v_n$ almost surely. Assume that $-1\leqslant a<1/2$ and let $\gamma>0$ be fixed. One can notice that 
\begin{equation}\label{equivalent_clef_a_plus_petit_que_un_demi}
	\frac{\sqrt{n}}{\left(\log n\right)^{\frac{1}{2}+\gamma}}\left(\frac{T_n}{n}-\frac{(1-\alpha)\mu_1}{1-a}\right)
	=\frac{\tilde{T}_n}{\sqrt{n}\left(\log n\right)^{\frac{1}{2}+\gamma}}=\frac{M_n}{a_n\sqrt{n}\left(\log n\right)^{\frac{1}{2}+\gamma}}.
\end{equation}
Moreover, using \eqref{estimee_de_v_n}, we have $\lim_{n\to+\infty}\langle M\rangle_n=\lim_{n\to+\infty}v_n=+\infty$. So, by \cite[Theorem 1.3.15]{Duflo1997}, we get 
\begin{equation}\label{convergence_M_n_a_plus_petit_que_un_demi}
\frac{M_n}{\langle M\rangle_n}=o\left(\sqrt{\frac{\left(\log \langle M\rangle_n\right)^{1+\gamma}}{\langle M\rangle_n}}\right)\,\textrm{a.s.}
\end{equation}
Combining \eqref{estimee_de_a_n}, \eqref{estimee_de_v_n}, \eqref{equivalent_clef_a_plus_petit_que_un_demi} and \eqref{convergence_M_n_a_plus_petit_que_un_demi}, we obtain 
$$
\frac{\sqrt{n}}{\left(\log n\right)^{\frac{1}{2}+\gamma}}\left(\frac{T_n}{n}-\frac{(1-\alpha)\mu_1}{1-a}\right)\converge{14}{n}{+\infty}{a.s.}0.
$$
Assume that $a=1/2$ and let $\gamma>0$ be fixed. Using \eqref{estimee_de_v_n}, we have $v_n\sim_{n\to+\infty}\frac{\pi}{4}\log n$. Moreover, since $M_n=a_n\tilde{T}_n$, we have 
\begin{equation}\label{equivalent_clef_a_egal_un_demi}
\frac{\sqrt{n}}{\sqrt{\log n}\left(\log \log n\right)^{\frac{1}{2}+\gamma}}\left(\frac{T_n}{n}-2(1-\alpha)\mu_1\right)=\frac{M_n}{a_n\sqrt{n\log n}\left(\log\log n\right)^{\frac{1}{2}+\gamma}}.
\end{equation}
Combining \eqref{estimee_de_a_n}, \eqref{estimee_de_v_n}, \eqref{convergence_M_n_a_plus_petit_que_un_demi} and \eqref{equivalent_clef_a_egal_un_demi}, we get
$$
\frac{\sqrt{n}}{\sqrt{\log n}\left(\log \log n\right)^{\frac{1}{2}+\gamma}}\left(\frac{T_n}{n}-2(1-\alpha)\mu_1\right)\converge{14}{n}{+\infty}{a.s.}0.
$$
Assume that $1/2<a<1$. Using \eqref{estimee_de_v_n}, we get $\lim_{n\to+\infty}v_n=\sum_{k=1}^{+\infty}a_k^2<+\infty$. Moreover, using \cite[Theorem 1.3.15]{Duflo1997}, we have $M_n\converge{10}{n}{+\infty}{a.s.}M_{\infty}:=\frac{(\alpha-a)\mu_1}{1-a}+\sum_{k=1}^{+\infty}a_k\varepsilon_k$. So, using \eqref{estimee_de_a_n}, we have 
$$
n^{1-a}\left(\frac{T_n}{n}-\frac{(1-\alpha)\mu_1}{1-a}\right)=\frac{\tilde{T}_n}{n^a}=\frac{M_n}{a_nn^a}\sim_{n\to+\infty}\frac{M_n}{\Gamma(a+1)}\converge{12}{n}{+\infty}{a.s.}L\quad\textrm{with}\quad L:=\frac{M_{\infty}}{\Gamma(a+1)}.
$$
Assume that $\E{\vert L\vert^\tau}<+\infty$ for some real $\tau\geqslant 2$ and let $n\geqslant 1$ be a fixed integer. Using Burkholder's inequality (\cite{Hall--Heyde1980}, Theorem 2.10), we get $\left\Vert M_n \right\Vert_\tau\lesssim\frac{(\alpha-a)\mu_1}{1-a}+\left(\sum_{k=1}^na_k^2\Vert\varepsilon_{k}\Vert_\tau^2\right)^{1/2}$. Keeping in mind that $\varepsilon_k=X_k-\mathbb E_{\F_{k-1}}[X_k]$ and using Lemma \ref{Lemme_technique}, we derive 
$$
\left\Vert M_n \right\Vert_\tau\lesssim\frac{(\alpha-a)\mu_1}{1-a}+\Vert\xi_1\Vert_\tau\sqrt{\sum_{k=1}^{+\infty}a_k^2}<+\infty.
$$ 
So, we obtain $\sup_{n\geqslant 1}\Vert M_n\Vert_\tau<+\infty$ and consequently  $\lim_{n\to+\infty}\mathbb E\left[\vert L_n-L\vert^{\tau}\right]=0$ where $L_n=\frac{M_n}{\Gamma(a+1)}$. Now, since $\mathbb E[\xi_1^2]$ is assumed to be finite, we are going to give an explicit formula for the first two moments of $L$. Since $(M_n)_{n\geqslant 1}$ is a martingale, we have 
\begin{equation}\label{moment_ordre_1_de_L}
\mathbb E[L]=\frac{\mathbb E[M_1]}{\Gamma(a+1)}=\frac{(\alpha-a) \mu_1}{(1-a) \Gamma(a+1)}.
\end{equation}
Since $\mathbb E_{\F_n}[\tilde{X}_{n+1}]=\frac{a}{n}\tilde{T}_n$, then it follows from the defintion of $\tilde{T}_n$ that  
$$
\mathbb E[\tilde{T}^2_{n+1}]=\left(1+\frac{2a}{n}\right)\mathbb E[\tilde{T}^2_n]+\mathbb E[\tilde{X}^2_{n+1}].
$$
Using Lemma $3.3$ in \cite{Kiss--Veto2022}, we derive that for any $n\geqslant 2$,
\begin{equation}\label{moment_ordre_deux_T_tilde_n}
\mathbb E[\tilde{T}^2_n]=\frac{\Gamma(n+2a)}{\Gamma(n)}\left(\frac{\mathbb E[\tilde{T}^2_1]}{\Gamma(1+2a)}+\sum_{j=1}^{n-1} \frac{\Gamma(j+1)\mathbb E[\tilde{X}^2_{j+1}]}{\Gamma(j+1+2a)}\right).
\end{equation}
By definition of $\tilde{X}_{n+1}$, we have 
\begin{equation}\label{esperance_X_n_plus_un}
\mathbb{E}[X_{n+1}]=\mathbb{E}[ \tilde{X}_{n+1}]+\frac{(1-\alpha) \mu_1}{1-a}=\frac{a}{n} \mathbb{E}[ \tilde{T}_n]+\frac{(1-\alpha) \mu_1}{1-a}.
\end{equation}
Moreover, since $\mathbb{E}_{\mathcal{F}_{n-1}}[\tilde{T}_n]=\left(1+\frac{a}{n}\right) \tilde{T}_n$, we have
\begin{equation}\label{esperance_T_tilde_n}
\mathbb{E}[\tilde{T}_n]=\mathbb{E}[\tilde{T}_1] \prod_{k=1}^{n-1}\left(1+\frac{a}{k}\right)=a_n^{-1} \mathbb{E}[\tilde{T}_1]=\frac{\Gamma(n+a) \mathbb{E}[\tilde{T}_1]}{\Gamma(n) \Gamma(a+1)}=\frac{(\alpha-a) \mu_1 \Gamma(n+a)}{(1-a) \Gamma(n) \Gamma(a+1)}.
\end{equation}
with the usual convention $\prod_{k=1}^0\left(1+\frac{a}{k}\right)=1$. Combining \eqref{esperance_X_n_plus_un} and \eqref{esperance_T_tilde_n}, we obtain
\begin{equation}\label{moment_ordre_deux_X_n_plus_un}
\mathbb E[X_{n+1}]=\frac{\mu_1}{1-a}\left(1-\alpha+\frac{(\alpha-a)\Gamma(n+a)}{\Gamma(a)\Gamma(n+1)}\right).
\end{equation}
Using Lemma \ref{Lemme_technique}, we have $\mathbb E[X^2_{n+1}]
=\mu_2$. Since 
$$
\mathbb{E}\left[\tilde{X}_{n+1}^2\right]=\mathbb{E}\left[X_{n+1}^2\right]-\frac{2(1-\alpha) \mu_1 \mathbb{E}\left[X_{n+1}\right]}{1-a}+\frac{(1-\alpha)^2 \mu_1^2}{(1-a)^2},
$$
we derive
$$
\mathbb{E}[\tilde{X}_{n+1}^2]=\sigma^2-\frac{c_0\Gamma(n+a)}{\Gamma(n+1)}
$$
where $\sigma^2$ is defined by \eqref{definition_sigma_carre_et_a} and 
\begin{equation}\label{definition_c_0}
c_0:=\frac{2(1-\alpha)(\alpha-a)\mu_1^2}{(1-a)^2\Gamma(a)}.
\end{equation}
Now, we have 
$$
\mathbb E[L_n^2]=\frac{a_n^2\mathbb E[\tilde{T}_n^2]}{\Gamma^2(a+1)}=\frac{\Gamma(n)\Gamma(n+2a)}{\Gamma^2(n+a)}\left(\frac{\mathbb E[\tilde{T}_1^2]}{\Gamma(1+2a)}+\sum_{j=1}^{n-1}\frac{\Gamma(j+1)\mathbb E[\tilde{X}_{j+1}^2]}{\Gamma(j+1+2a)}\right)
$$
where 
$$
\mathbb E[\tilde{T}^2_1]
=\mathbb E\left[\left(X_1-\frac{(1-\alpha)\mu_1}{1-a}\right)^2\right]=\mu_2-\frac{(1-\alpha)(1+\alpha-2a)\mu_1^2}{(1-a)^2}.
$$
Consequently,
\begin{align*}
	\frac{\Gamma^2(n+a)\mathbb E[L_n^2]}{\Gamma(n)\Gamma(n+2a)}&=
	\frac{1}{\Gamma(2a+1)}\left(\mu_2-\frac{(1-\alpha)(1+\alpha-2a)\mu_1^2}{(1-a)^2}\right)\\
	&\qquad+\sigma^2\sum_{j=1}^{n-1}\frac{\Gamma(j+1)}{\Gamma(j+1+2a)}-c_0\sum_{j=1}^{n-1}\frac{\Gamma(j+a)}{\Gamma(j+1+2a)}.
\end{align*}
\begin{Lem} (Bercu \cite{Bercu2018})\label{Lemme_somme_de_Gamma} {\em For any non-negative real numbers $s$ and $t$ such that $t \neq s+1$ and for all $n \geq 1$, we have}
	$$
	\sum_{k=1}^n \frac{\Gamma(k+s)}{\Gamma(k+t)}=\frac{\Gamma(n+s+1)}{(t-s-1) \Gamma(n+t)}\left(\frac{\Gamma(n+t) \Gamma(s+1)}{\Gamma(n+s+1) \Gamma(t)}-1\right) .
	$$
\end{Lem}
\noindent Since $1/2<a<1$, using Lemma \ref{Lemme_somme_de_Gamma}, for $n\geqslant 2$, we have 
\begin{equation}\label{somme1_Gamma}
\sum_{j=1}^{n-1}\frac{\Gamma(j+1)}{\Gamma(j+1+2a)}=\frac{\Gamma(n+1)}{(2a-1)\Gamma(n+2a)}\left(\frac{\Gamma(n+2a)}{\Gamma(n+1)\Gamma(2a+1)}-1\right),
\end{equation}
and
\begin{equation}\label{somme2_Gamma}
\sum_{j=1}^{n-1}\frac{\Gamma(j+a)}{\Gamma(j+1+2a)}=\frac{\Gamma(n+a)}{a\Gamma(n+2a)}\left(\frac{\Gamma(n+2a)\Gamma(a+1)}{\Gamma(n+a)\Gamma(1+2a)}-1\right).
\end{equation}	
Using \eqref{somme1_Gamma} and \eqref{somme2_Gamma}, we obtain 
\begin{align*}
\frac{\Gamma^2(n+a)\mathbb E[L_n^2]}{\Gamma(n)\Gamma(n+2a)}
&=\frac{1}{\Gamma(2a+1)}\left(\mu_2-\frac{(1-\alpha)(1+\alpha-2a)\mu_1^2}{(1-a)^2}\right)\\
&\qquad+\frac{\sigma^2\Gamma(n+1)}{(2a-1)\Gamma(n+2a)}\left(\frac{\Gamma(n+2a)}{\Gamma(n+1)\Gamma(2a+1)}-1\right)\\
&\quad\qquad\qquad-\frac{c_0\Gamma(n+a)}{a\Gamma(n+2a)}\left(\frac{\Gamma(n+2a)\Gamma(a+1)}{\Gamma(n+a)\Gamma(2a+1)}-1\right).
\end{align*}
Since $\mathbb E[L^2]=\lim_{n\to+\infty}\mathbb E[L_n^2]$, we derive 
\begin{equation}\label{moment_ordre_2_de_L}
\mathbb E[L^2]
=\frac{1}{\Gamma(2a+1)}\left(\mu_2-\frac{(1-\alpha)(1+\alpha-2a)\mu_1^2}{(1-a)^2}\right)+\frac{\sigma^2}{(2a-1)\Gamma(2a+1)}-\frac{2(1-\alpha)(\alpha-a)\mu_1^2}{(1-a)^2\Gamma(2a+1)}.
\end{equation}
Now, we are going to derive necessary and sufficient conditions for $L$ to be non-degenerate. In fact, we have 
$$
L=\frac{\mathbb{E} M_1}{\Gamma(a+1)}+\frac{1}{\Gamma(a+1)} \sum_{k=1}^{\infty} a_k \varepsilon_k
$$
where $\varepsilon_1=X_1-\mu_1$ and 
$
\varepsilon_k=X_k-\mathbb{E}_{\F_{k-1}}[X_k ]=X_k-\frac{a }{k-1}\sum_{j=1}^{k-1} X_j-(1-\alpha) \mu_1
$
for any $k\geqslant 2$. Thus, $L$ is degenerate if and only if $\varepsilon_k=0$ a.s. for any $k\geqslant 1$. If $\mu_2>\mu_1^2$, then $\varepsilon_1$ is non-degenerate, and thus $L$ is non-degenerate. Now assume that $\mu_2=\mu_1^2$, then $\xi_n=\mu_1$ a.s. for all $n \geq 1$. Note that $X_k=0$ a.s. for any $k\geqslant 1$ if $\mu_1=0$, in which case $L$ is degenerate. By \eqref{def_sympa_step-reinforced_random_walk},
$$
\varepsilon_2=\gamma_1 \rho_1 \mu_1+\left(1-\gamma_1\right) \mu_1-a \mu_1-(1-\alpha) \mu_1=\mu_1\left(\gamma_1(\rho_1-1)+\alpha-a\right),
$$
where $\gamma_1$ is a Bernoulli random variable with parameter $\alpha$ and $\rho_1$ is a Rademacher random variable with parameter $p$. Moreover, $\gamma_1$ and $\rho_1$ are independent. If $\mu_1 \neq 0$, then $\varepsilon_2=0$ a.s. if and only if
$
\left(\rho_1-1\right) \gamma_1=a-\alpha>-1/2
$ a.s. which occurs if and only if $p=1$ since otherwise the left-hand side equals $-2$ with positive probability. In this case (i.e. $p=1$), $X_k=\mu_1$ a.s. for all $k \geq 1$, and in particular, $L$ is degenerate. In summary, $L$ is degenerate if and only if $\mu_2=\mu_1^2$ and either $\mu_1=0$ or $p=1$. The proof of Theorem \ref{Loi_des_grands_nombres_sous_moment_ordre_deux} is complete.$\hfill\Box$\\
\\
{\em Proof of Theorem \ref{TLC}}. Recall that $M_n=\frac{(\alpha-a)\mu_1}{1-a}+\tilde{M}_n+\overline{M}_n$ with $\tilde{M}_n=\sum_{k=1}^na_k\tilde{\varepsilon}_k$ and $\overline{M}_n=\sum_{k=1}^na_k\overline{\varepsilon}_k$ where $\tilde{\varepsilon}_k$ and $\overline{\varepsilon}_k$ are defined by \eqref{def_epsilon_barre_k} from an increasing sequence $(c_k)_{k\geqslant 1}$ of positive real numbers satisfying $\lim_{k\to+\infty}c_k=+\infty$. Assume that $-1\leqslant a<1/2$. Using \eqref{estimee_de_v_n}, we have  $v_n:=\sum_{k=1}^na_k^2\to+\infty$ as $n\to+\infty$. Using Lemma \ref{Convergence_de_moment_ordre_2_de_epsilon_barre_k} and noting that $\mathbb E[\overline{M}_n^2]=\sum_{k=1}^na_k^2\mathbb E[\overline{\varepsilon}_k^2]$, we get $\lim_{n\to+\infty}v_n^{-1}\mathbb E[ \overline{M}_n^2]=0$. In particular, we obtain $v_n^{-1/2}\overline{M}_n\converge{10}{n}{+\infty}{\textrm{$\mathbb P$}}0$. Moreover, from Lemma \ref{Convergence_de_moment_ordre_2_de_epsilon_barre_k}, we have $\mathbb E_{\F_{k-1}}\left[\tilde{\varepsilon}_k^2\right]\converge{10}{k}{+\infty}{\textrm{$\mathbb L^1$}}\sigma^2$ and consequently $
v_n^{-1}\sum_{k=1}^na_k^2\,\mathbb E_{\F_{k-1}}[\tilde{\varepsilon}_k^2]\converge{12}{n}{+\infty}{\textrm{$\mathbb L^1$}}\sigma^2$ where $\sigma^2$ is given by \eqref{definition_sigma_carre_et_a}. From now on, we assume that $c_k=k^{\rho}$ with $0<\rho<1/2$ for any integer $k\geqslant 1$. Let $\eta>0$ be fixed and recall that $\vert\tilde{\varepsilon}_k\vert\leqslant 2c_k$ a.s. and $\mathbb E[\tilde{\varepsilon}_k^2]\leqslant \mathbb E[\xi_1^2]<+\infty$. Then, using \eqref{estimee_de_a_n} and \eqref{estimee_de_v_n}, we get
$$
	\frac{1}{v_n}\sum_{k=1}^n\mathbb E[a_k^2\tilde{\varepsilon}_k^2\ind_{\vert a_k\tilde{\varepsilon}_k\vert>\eta\sqrt{v_n}}]\leqslant\frac{2\mathbb E[\xi_1^2]}{\eta v_n^{3/2}}\sum_{k=1}^na_k^3c_k
	\lesssim \left\{\begin{array}{ccc}
		n^{-\frac{3}{2}+3a}\sum_{k=1}^{+\infty}k^{-3a+\rho} & \textrm{if} & 3a-\rho>1\\
 	n^{-\frac{3}{2}+3a}\log n & \textrm{if} & 3a-\rho=1\\
 		n^{-\frac{1}{2}+\rho} & \textrm{if} & 3a-\rho<1
	\end{array}\right.
$$
Consequently,
\begin{equation}\label{condition_de_Lindeberg}
\frac{1}{v_n}\sum_{k=1}^n\mathbb E_{\F_{k-1}}[a_k^2\tilde{\varepsilon}_k^2\ind_{\vert a_k\tilde{\varepsilon}_k\vert>\eta\sqrt{v_n}}]\converge{15}{n}{+\infty}{$\mathbb L^1$}0.
\end{equation}	
Applying Corollary 2.1.10 in \cite{Duflo1997}, we derive 
$$
\frac{\tilde{M}_n}{\sqrt{v_n}}\converge{12}{n}{+\infty}{Law}\mathcal N\left(0,\sigma^2\right).
$$
Using Slutsky's lemma, we obtain $v_n^{-1/2}M_n\converge{10}{n}{+\infty}{\textrm{Law}}\mathcal N\left(
0,\sigma^2\right)$. Since $M_n=a_n\tilde{T}_n$, it means that 
$$
\dfrac{na_n}{\sqrt{v_n}} \left(\frac{T_n}{n}-\frac{(1-\alpha)\mu_1}{1-a}\right)\converge{12}{n}{+\infty}{\textrm{Law}}\mathcal N\left(
0,\sigma^2\right).
$$
Using \eqref{estimee_de_a_n} and \eqref{estimee_de_v_n}, we obtain   
$$
\sqrt{n} \left(\frac{T_n}{n} - \frac{(1-\alpha)\mu_1}{1-a} \right) \converge{12}{n}{+\infty}{\textrm{Law}}\mathcal{N} \left(0, \frac{\sigma^2}{1-2a} \right).
$$
Assume that $a=1/2$. From Lemma \ref{Convergence_de_moment_ordre_2_de_epsilon_barre_k}, we have $\lim_{k\to+\infty}\mathbb E[\overline{\varepsilon}_k^2]=0$ and using \eqref{estimee_de_v_n}, we have $v_n\sim_{n\to+\infty}\frac{\pi}{4}\log n$. So,
$$
v_n^{-1}\mathbb E[\overline{M}_n^2]=v_n^{-1}\sum_{k=1}^na_k^2\mathbb E[\overline{\varepsilon}_k^2]\converge{12}{n}{+\infty}{}0.
$$
From Lemma \ref{Convergence_de_moment_ordre_2_de_epsilon_barre_k}, we have also $\mathbb E_{\F_{k-1}}[\tilde{\varepsilon}^2_k]\converge{10}{k}{+\infty}{\textrm{$\mathbb L^1$}}\sigma^2$ and consequently, we derive 
$$
\frac{1}{v_n}\sum_{k=1}^na_k^2\mathbb E_{\F_{k-1}}[\tilde{\varepsilon}_k^2]\converge{12}{n}{+\infty}{\textrm{$\mathbb L^1$}}\sigma^2.
$$
As before, we assume that $c_k=k^{\rho}$ with $0<\rho<1/2$ for any integer $k\geqslant 1$. If $\eta>0$ is fixed, then, using \eqref{estimee_de_a_n} and \eqref{estimee_de_v_n}, we have 
$$
\frac{1}{v_n}\sum_{k=1}^n\mathbb E[a_k^2\tilde{\varepsilon}_k^2\ind_{\vert a_k\tilde{\varepsilon}_k\vert>\eta\sqrt{v_n}}]\leqslant\frac{2\mathbb E[\xi_1^2]}{\eta v_n^{3/2}}\sum_{k=1}^na_k^3c_k
\lesssim 
	(\log n)^{-3/2}\sum_{k=1}^{+\infty}k^{-\frac{3}{2}+\rho}\converge{12}{n}{+\infty}{}0. 
$$
So, we get 
$$
\frac{1}{v_n}\sum_{k=1}^n\mathbb E_{\F_{k-1}}[a_k^2\tilde{\varepsilon}_k^2\ind_{\vert a_k\tilde{\varepsilon}_k\vert>\eta\sqrt{v_n}}]\converge{15}{n}{+\infty}{$\mathbb L^1$}0.
$$
Using Corollary 2.1.10 in \cite{Duflo1997}, we derive $v_n^{-1/2}\tilde{M}_n\converge{10}{n}{+\infty}{\textrm{Law}}\mathcal N(0,\sigma^2)$. Applying Slutsky's lemma, we obtain $v_n^{-1/2}M_n\converge{10}{n}{+\infty}{\textrm{Law}}\mathcal N(0,\sigma^2)$ which means that
$$\dfrac{na_n}{\sqrt{v_n}} \times \left(\frac{T_n}{n}-2(1-\alpha)\mu_1 \right)\converge{10}{n}{+\infty}{\textrm{Law}}\mathcal N\left(
0,\sigma^2\right).
$$
So, using \eqref{estimee_de_a_n} et \eqref{estimee_de_v_n}, we obtain   
$$
\frac{\sqrt{n}}{\sqrt{\log{n}}} \left(\frac{T_n}{n} - 2(1-\alpha)\mu_1 \right) \converge{12}{n}{+\infty}{\textrm{Law}}\mathcal{N} \left(0, \sigma^2 \right)\quad\textrm{with}\quad \sigma^2=\mu_2-4(1-\alpha)^2\mu_1^2.
$$
Assume that $1/2<a<1$. Using \eqref{estimee_de_v_n}, we derive $\lim_{n\to+\infty}\sum_{k=n+1}^{+\infty}a_k^2=0$. First, we are going to prove  
\begin{equation}\label{TLC_a_la_Heyde}
	\frac{\sum_{k=n+1}^{+\infty}a_k\varepsilon_k}{\sqrt{\sum_{k=n+1}^{+\infty}a_k^2\mathbb E[\varepsilon_k^2]}}\converge{12}{n}{+\infty}{\textrm{Law}}\mathcal N\left(0,1\right).	
\end{equation}
Using Lemma \ref{Lemme_technique} and Lemma \ref{Convergence_de_moment_ordre_2_de_epsilon_barre_k}, we have 
$$
\mathbb E[\varepsilon_k^2]=\mathbb E[\xi_1^2]-\mathbb E\!\left[\left(\mathbb E_{\F_{k-1}}[X_k]\right)^2\right]\converge{12}{k}{+\infty}{}\mu_2-\frac{(1-\alpha)^2\mu_1^2}{(1-a)^2}=\sigma^2.
$$ 
From Lemma \ref{Convergence_de_moment_ordre_2_de_epsilon_barre_k}, we have also $\lim_{k\to+\infty}\mathbb E[\overline{\varepsilon}_k^2]=0$. So, we get 
$$
\left\Vert\frac{\sum_{k=n+1}^{+\infty}a_k\overline{\varepsilon}_k}{\sqrt{\sum_{k=n+1}^{+\infty}a_k^2\mathbb E[\varepsilon_k^2]}}\right\Vert_2^2=\frac{\sum_{k=n+1}^{+\infty}a_k^2\mathbb E[\overline{\varepsilon}_k^2]}{\sum_{k=n+1}^{+\infty}a_k^2\mathbb E[\varepsilon_k^2]}\converge{12}{n}{+\infty}{}0.
$$
So, in order to get \eqref{TLC_a_la_Heyde},  it is sufficient to prove 
\begin{equation}\label{TLC_a_la_Heyde_tilde}
	\frac{\sum_{k=n+1}^{+\infty}a_k\tilde{\varepsilon}_k}{\sqrt{\sum_{k=n+1}^{+\infty}a_k^2}}\converge{12}{n}{+\infty}{\textrm{Law}}\mathcal N\left(0,\sigma^2\right).	
\end{equation}
First, 
$$
\left\Vert\frac{\sum_{k=n+1}^{+\infty}a_k^2\left(\tilde{\varepsilon}_k^2-\mathbb E_{\F_{k-1}}[\tilde{\varepsilon}_k^2]\right)}{\sum_{k=n+1}^{+\infty}a_k^2}\right\Vert_2^2=s_{n+1}^{-4}\sum_{k=n+1}^{+\infty}a_k^4\mathbb E\left[\left(\tilde{\varepsilon}_k^2-\mathbb E_{\F_{k-1}}[\tilde{\varepsilon}_k^2]\right)^2\right]
$$
where $s_n^2:=\sum_{k=n}^{+\infty}a_k^2$ satisfies $\lim_{n\to+\infty}(2a-1)n^{2a-1}s_n^2=\Gamma^2(1+a)$. Again, we assume that $c_k=k^{\rho}$ with $0<\rho<1/2$ for any integer $k\geqslant 1$. Using \eqref{estimee_de_a_n} and \eqref{estimee_de_v_n} and keeping in mind that  $\vert\tilde{\varepsilon}_k\vert\leqslant 2c_k$ a.s. and $\mathbb E[\tilde{\varepsilon}_k^2]\leqslant \mathbb E[\xi_1^2]<+\infty$, we get
\begin{equation}\label{convergence_L2_epsilon_tilde_carre}
	\left\Vert\frac{\sum_{k=n+1}^{+\infty}a_k^2\left(\tilde{\varepsilon}_k^2-\mathbb E_{\F_{k-1}}[\tilde{\varepsilon}_k^2]\right)}{\sum_{k=n+1}^{+\infty}a_k^2}\right\Vert_2^2
	\lesssim \frac{\mathbb E[\xi_1^2]}{n^{2-4a}}\sum_{k=n+1}^{+\infty}k^{-4a+2\rho}
	\leqslant\mathbb E[\xi_1^2] \sum_{k=n+1}^{+\infty}\left(\frac{1}{k}\right)^{2(1-\rho)}\converge{10}{n}{+\infty}{}0.
\end{equation}
From Lemma \ref{Convergence_de_moment_ordre_2_de_epsilon_barre_k}, we have $\mathbb E_{\F_{k-1}}[\tilde{\varepsilon}_k^2]\converge{10}{k}{+\infty}{\textrm{$\mathbb L^1$}}\sigma^2$ and consequently 
\begin{equation}\label{convergence_L2_epsilon_tilde_carre_2eme_partie}
	s_n^{-2}\sum_{k=n}^{+\infty}a_k^2\left(\mathbb E_{\F_{k-1}}[\tilde{\varepsilon}_k^2]-\sigma^2\right)\converge{12}{n}{+\infty}{\textrm{$\mathbb L^1$}}0.	
\end{equation}
Combining \eqref{convergence_L2_epsilon_tilde_carre} and \eqref{convergence_L2_epsilon_tilde_carre_2eme_partie}, we derive 
\begin{equation}\label{convergence_L1_moyennes_des_epsilon_tilde_carre_vers_sigma_carre}
	\frac{\sum_{k=n}^{+\infty}a_k^2\tilde{\varepsilon}_k^2}{\sum_{k=n}^{+\infty}a_k^2}\converge{14}{n}{+\infty}{\textrm{$\mathbb L^1$}}\sigma^2.	
\end{equation}
Now, we are going to prove the following Lindeberg's type condition 
\begin{equation}\label{convergence_de_I_n}
	I_n:=\frac{1}{s_n^2}\sum_{k=n}^{+\infty}\mathbb E[a_k^2\tilde{\varepsilon}_k^2\ind_{\vert a_k\tilde{\varepsilon}_k\vert>\delta s_n}]\converge{12}{n}{+\infty}{}0.
\end{equation}
for any $\delta>0$. Using again \eqref{estimee_de_a_n} and \eqref{estimee_de_v_n}, we get
$$
I_n\leqslant \frac{2\mathbb E[\xi_1^2]}{\delta s_n^{3}}\sum_{k=n}^{+\infty}a_k^3c_k\lesssim n^{3a-\frac{3}{2}}\sum_{k=n}^{+\infty}k^{-3a+\rho} \leqslant \sum_{k=n}^{+\infty}k^{-\frac{3}{2}+\rho}.
$$
Since $3/2-\rho>1$, we obtain \eqref{convergence_de_I_n}. So, using Theorem \textrm{A$2$} in \cite{Kubota--Takei2019}, we obtain \eqref{TLC_a_la_Heyde_tilde}. Using Slutsky's lemma, we derive 
\begin{equation}\label{convergence_en_loi_finale}
	\frac{\sum_{k=n+1}^{+\infty}a_k\varepsilon_k}{\sqrt{\sum_{k=n+1}^{+\infty}a_k^2}}\converge{12}{n}{+\infty}{\textrm{Law}}\mathcal N\left(0,\sigma^2\right).	
\end{equation}
Since $\lim_{n\to+\infty}(2a-1)n^{2a-1}s_n^2=\Gamma^2(1+a)$, one can notice that \eqref{convergence_en_loi_finale} means that 
\begin{equation}\label{convergence_en_loi_finale2}
	\sqrt{n^{2a-1}}\left(\frac{n a_n}{\Gamma(1+a)}\left(\frac{T_n}{n}-\frac{(1-\alpha)\mu_1}{1-a}\right)-L\right)\converge{12}{n}{+\infty}{\textrm{Law}}\mathcal N\left(0,\frac{\sigma^2}{2a-1}\right)
\end{equation}
where $L=\frac{M_{\infty}}{\Gamma(a+1)}$  and  $M_{\infty}=\frac{(\alpha-a)\mu_1}{1-a}+\sum_{k=1}^{+\infty}a_k\varepsilon_k$. Finally, using  \eqref{estimee_de_v_n} and \eqref{convergence_en_loi_finale2} and 
$$
\frac{a_n}{\Gamma(a+1)}=\frac{1}{n^a}\left(1+\frac{a(1-a)}{2 n}+O\left(\frac{1}{n^2}\right)\right),
$$ 
we obtain
$$
\sqrt{n^{2a-1}}\left(n^{1-a}\left(\frac{T_n}{n}-\frac{(1-\alpha)\mu_1}{1-a}\right)-L\right)\converge{12}{n}{+\infty}{\textrm{Law}}\mathcal N\left(0,\frac{\sigma^2}{2a-1}\right)
$$ 
The proof of Theorem \ref{TLC} is complete.$\hfill\Box$\\
\\
{\em Proof of Theorem \ref{TLC_bis}}. Let $-1\leqslant a<1$ be fixed. Let also $n\geqslant 1$ be a fixed integer and denote $\tau_n:=n^{-1}m_n+a\left(1-n^{-1}m_n\right)$ and $\kappa_n:=n^{-1}m_n(1-n^{-1}m_n)$. For any integer $k\geqslant 1$, we consider the $\sigma$-algebras $\F_k:=\sigma(X_j\,;\,1\leqslant j\leqslant k)$ and
$\F'_{k,n}:=\F_{k}\vee \sigma\left(X'_j\,;\,m_n<j\leqslant k\right)
$ with the convention $\F'_{k,n}=\F_k$ if $1\leqslant k\leqslant m_n$. For simplicity, we write $\F'_k$ instead of $\F'_{k,n}$. Using \eqref{definition_unbalanced_step-reinforced_random_walk} and \eqref{Def_X_prime}, for any $k\geqslant 2$, we have  
$\mathbb E_{\F'_{k-1}}[X_k]=\mathbb E_{\F_{k-1}}[X_k]=a(k-1)^{-1}T_{k-1}+(1-\alpha)\mu_1$ and, for any $k>m_n$, $\mathbb E_{\F'_{k-1}}[X'_k]=\mathbb E_{\F_{m_n}}[X'_k]=am_n^{-1}T_{m_n}+(1-\alpha)\mu_1$. Then,
$$
\sqrt{m_n}\,\frac{ S_n}{n}=\tau_n \sqrt{m_n}\left(\frac{T_{m_n}}{m_n}- \frac{(1-\alpha)\mu_1}{1-a}\right) +\frac{\sqrt{m_n}}{n} \sum_{k=m_n+1}^n\left(X'_k-\mathbb E_{\F'_{k-1}}[X'_k]\right)+\frac{ \sqrt{m_n}(1-\alpha)\mu_1}{1-a}
$$
and
$$
\sqrt{m_n}\left(\frac{ S_n}{n}-\frac{(1-\alpha)\mu_1}{1-a}\right)=G_n+\sqrt{\kappa_n}H_n
$$
where 
\begin{equation}\label{def_G_n_et_H_n}
G_n=\tau_n \sqrt{m_n}\left(\frac{T_{m_n}}{m_n}- \frac{(1-\alpha)\mu_1}{1-a}\right) \quad\textrm{and}\quad H_n=\frac{1}{\sqrt{n-m_n}} \sum_{k=m_n+1}^nY'_k
\end{equation}
with $Y'_k=X'_k-\mathbb E_{\F'_{k-1}}[X'_k]=X'_k-\mathbb E_{\F_{m_n}}[X'_k]$ for any $m_n<k\leqslant n$. Now, we consider the charactersitic function $\psi_n$ defined for any $t\in\R$ by $\psi_n(t)=\mathbb E\left[\exp\left(it\left(G_n+\sqrt{\kappa_n}H_n\right)\right)\right]$. Let $\ell\geqslant 2$ and $i_{\ell}>i_{\ell-1}>\dots>i_1>i_0=m_n$ be fixed. For any $1\leqslant s\leqslant \ell$, we have $\mathbb E_{\F'_{i_{s-1}}}[\varphi_s(X'_{i_{s}})]=\mathbb E_{\F_{m_n}}[\varphi_s(X'_{i_{s}})]$ for any measurable and bounded function $\varphi_s$. Consequently,
\begin{align*}
\mathbb E_{\F_{m_n}}\left[\prod_{s=1}^{\ell}\varphi_s\left(X'_{i_{s}}\right)\right]&=\mathbb E_{\F_{m_n}}\left[\mathbb E_{\F'_{i_{\ell-1}}}[\varphi_{\ell}\left(X'_{i_{\ell}}\right)]\prod_{s=1}^{\ell-1}\varphi_s\left(X'_{i_{s}}\right)\right]=\mathbb E_{\F_{m_n}}[\varphi_{\ell}\left(X'_{i_{\ell}}\right)]\times\mathbb E_{\F_{m_n}}\left[\prod_{s=1}^{\ell-1}\varphi_s\left(X'_{i_{s}}\right)\right]
\end{align*}
and, by induction, we have
$$
\mathbb E_{\F_{m_n}}\left[\prod_{s=1}^{\ell}\varphi_s\left(X'_{i_{s}}\right)\right]=\prod_{s=1}^{\ell}\mathbb E_{\F_{m_n}}[\varphi_s\left(X'_{i_{s}}\right)].
$$
So, conditionally to the $\sigma$-algebra $\F_{m_n}$, the random variables $(X'_k)_{m_n<k\leqslant n}$ are independent and identically distributed. Moreover, for any $m_n<k\leqslant n$, we have $\mathbb E_{\F_{m_n}}[Y'_k]=0$ and 
\begin{align*}
\mathbb V_{\F_{m_n}}[Y'_k]&=\mathbb E_{\F_{m_n}}[X'^2_k]-\left(\mathbb E_{\F_{m_n}}[X'_k]\right)^2=\frac{\alpha}{m_n}\sum_{\ell=1}^{m_n}X_{\ell}^2+(1-\alpha)\mu_2-\left(\frac{a}{m_n}\sum_{\ell=1}^{m_n}X_{\ell}+(1-\alpha)\mu_1\right)^2.
\end{align*}
In fact, we have $\mathbb V_{\F_{m_n}}[Y'_k]=\mathbb E_{\F_{m_n}}[\varepsilon^2_{m_n+1}]$ for any $m_n<k\leqslant n$ where $\varepsilon_{\ell}=X_{\ell}-\mathbb E_{\F_{\ell-1}}[X_{\ell}]$ for any $\ell\geqslant 1$. Using Lemma \ref{Convergence_de_moment_ordre_2_de_epsilon_barre_k}, we have
\begin{equation}\label{convergence_variance_conditionnelle_de_Y_prime}
\mathbb V_{\F_{m_n}}[Y'_{m_n+1}]\converge{19}{n}{+\infty}{$\mathbb L^1$ and a.s.}\sigma^2
\end{equation}
where $\sigma^2$ is given by \eqref{definition_sigma_carre_et_a}. Now, we note that $\psi_n(t)=\mathbb E[\exp\left(itG_n\right)\mathbb E_{\F_{m_n}}[\exp\left(it\sqrt{\kappa_n}H_n\right)]]$ for any $t\in\R$. In order to derive the limit of $\psi_n(t)$ as $n$ goes to infinity, we have to prove that
\begin{equation}\label{Lemme_clef_convergence_fonction_caracteristique_conditionnelle}
\mathbb E_{\F_{m_n}}[\exp\left(it\sqrt{\kappa_n}H_n\right)]\converge{19}{n}{+\infty}{$\mathbb L^1$ and a.s.}\exp\left(-\frac{\theta(1-\theta)\sigma^2t^2}{2}\right).
\end{equation}
In fact, if $n\geqslant 1$ is a fixed integer, then
\begin{align*}
&\mathbb E_{\F_{m_n}}[\exp\left(it\sqrt{\kappa_n}H_n\right)]=\prod_{k=m_n+1}^n\mathbb E_{\F_{m_n}}\left[\exp\left(\frac{it\sqrt{\kappa_n}Y'_k}{\sqrt{n-m_n}}\right)\right]\\
&=\prod_{k=m_n+1}^n\mathbb E_{\F_{m_n}}\left[1+\frac{it\sqrt{\kappa_n}Y'_k}{\sqrt{n-m_n}}-\frac{t^2\kappa_nY'^2_k}{2(n-m_n)}+R_n\right]=\left(1-\frac{t^2\kappa_n\mathbb E_{\F_{m_n}}[Y'^2_{m_n+1}]}{2(n-m_n)}+\mathbb E_{\F_{m_n}}[R_n]\right)^{n-m_n}
\end{align*}
where, using \cite[Lemma 1.2, page 556]{Gut2013}, we have
$$
\vert R_n\vert\leqslant \frac{t^2\kappa_nY'^2_{m_n+1}}{n-m_n}\left(1\wedge \frac{t\vert Y'_{m_n+1}\vert\sqrt{\kappa_n}}{\sqrt{n-m_n}}\right).
$$
Moreover, using $Y'^2_{m_n+1}\leqslant 2(X'^2_{m_n+1}+\left(\mathbb E_{\F_{m_n}}[X'_{m_n+1}]\right)^2)$, we have 
$$
(n-m_n)\mathbb E_{\F_{m_n}}[\vert R_n\vert]
\lesssim\mathbb E_{\F_{m_n}}\left[ Y'^2_{m_n+1}\left(1\wedge \frac{\vert Y'_{m_n+1}\vert}{\sqrt{n-m_n}}\right)\right]\lesssim A_{1,n}+A_{2,n}+A_{3,n}
$$
where 
\begin{align*}
A_{1,n}&=\mathbb E_{\F_{m_n}}\left[ X'^2_{m_n+1}\left(1\wedge \frac{\vert X'_{m_n+1}\vert}{\sqrt{n-m_n}}\right)\right]\\
A_{2,n}&=\mathbb E_{\F_{m_n}}\left[ X'^2_{m_n+1}\left(1\wedge \frac{\mathbb E_{\F_{m_n}}[\vert X'_{m_n+1}\vert]}{\sqrt{n-m_n}}\right)\right]\\
A_{3,n}&=\mathbb E_{\F_{m_n}}\left[\left(\mathbb E_{\F_{m_n}}[ X'_{m_n+1}]\right)^2\left(1\wedge \frac{\vert Y'_{m_n+1}\vert}{\sqrt{n-m_n}}\right)\right].
\end{align*}
Using Lemma \ref{Convergence_de_moment_ordre_2_de_epsilon_barre_k}, we have 
$$
\mathbb E_{\F_{m_n}}[X'^2_{m_n+1}]=\frac{\alpha}{m_n}\sum_{\ell=1}^{m_n}X_{\ell}^2+(1-\alpha)\mu_2\converge{12}{n}{+\infty}{a.s.}\mu_2.
$$
Moreover, using \eqref{convergence_variance_conditionnelle_de_Y_prime}, we have also $\mathbb E_{\F_{m_n}}[Y'^2_{m_n+1}]=\mathbb V_{\F_{m_n}}[Y'_{m_n+1}]\converge{12}{n}{+\infty}{a.s.}\sigma^2$. So, we get  
$$
\vert A_{3,n}\vert\leqslant \mathbb E_{\F_{m_n}}[X'^2_{m_n+1}]\frac{\sqrt{\mathbb E_{\F_{m_n}}[Y'^2_{m_n+1}]}}{\sqrt{n-m_n}}\converge{15}{n}{+\infty}{a.s.}0
$$
and
$$
\vert A_{2,n}\vert\leqslant \mathbb E_{\F_{m_n}}[
X'^2_{m_n+1}]\left(1\wedge \frac{\sqrt{\mathbb E_{\F_{m_n}}[ X'^2_{m_n+1}]}}{\sqrt{n-m_n}}\right)\converge{15}{n}{+\infty}{a.s.}0.
$$
Let $\phi_n$ be the function defined for any $x\in\R$ by $\phi_n(x)=x^2\left(1\wedge\frac{\vert x\vert}{\sqrt{n-m_n}}\right)$. Keeping in mind \eqref{def_sympa_step-reinforced_random_walk}, we have
\begin{align*}
A_{1,n}&=\mathbb E_{\F_{m_n}}[\phi_n(X'_{m_n+1})]=\mathbb E_{\F_{m_n}}[\phi_n\left(X_{U_{m_n+1,n}}\right)\ind_{\{\gamma_{m_n+1}=1\}}]+
\mathbb E_{\F_{m_n}}[\phi_n\left(\xi_{m_n+1}\right)\ind_{\{\gamma_{m_n+1}=0\}}]\\
&=\sum_{\ell=1}^{m_n}\mathbb E_{\F_{m_n}}[\phi_n(X_{\ell})\ind_{\{U_{m_n+1,n}=\ell\}}\ind_{\{\gamma_{m_n+1}=1\}}]
+(1-\alpha)\mathbb E[\phi_n(\xi_1)]\\
&=(1-\alpha)\mathbb E[\phi_n(\xi_1)]+\frac{\alpha}{m_n}\sum_{\ell=1}^{m_n}\phi_n(X_{\ell})\\
&=\left(1-\alpha\right)\mathbb E\left[\xi_1^2\left(1\wedge \frac{\vert\xi_1\vert}{\sqrt{n-m_n}}\right)\right]+\frac{\alpha}{m_n}\sum_{\ell=1}^{m_n}X_{\ell}^2\left(1\wedge \frac{\vert X_{\ell}\vert}{\sqrt{n-m_n}}\right).
\end{align*}
From the Lebesgue's dominated convergence theorem, we have 
$$
\lim_{n\to+\infty}\mathbb E\left[\xi_1^2\left(1\wedge \frac{\vert\xi_1\vert}{\sqrt{n-m_n}}\right)\right]=0.
$$
If $D_k:=\{|X_k|>\sqrt{k}\}$ for any $k\geqslant 1$, then 
$
\sum_{k=1}^{+\infty}\P(D_k)=\sum_{k=1}^{+\infty}\P\left(\xi_1^2> k\right)\leqslant\mathbb E\left[\xi_1^2\right]<+\infty$. Using Borel-Cantelli's lemma, we obtain $\P\left(N\right)=0$ where $N:=\cap_{n\geq1}\cup_{k\geq n}D_k$. Let $\omega\notin N$ be fixed. There exists a positive integer $K=K_\omega$ such that for any $k\geqslant K$, we have $\vert X_k(\omega)\vert\leqslant \sqrt{k}$. It means that for almost all $\omega\in\Omega$ and any $n\geqslant 1$ such that $m_n\geqslant K$,
\begin{equation}\label{control_A_1_n_part1}
\frac{1}{m_n}\sum_{\ell=1}^{m_n}X_{\ell}^2(\omega)\ind_{\{\vert X_{\ell}(\omega)\vert> \sqrt{\ell}\}}\left(1\wedge \frac{\vert X_{\ell}(\omega)\vert}{\sqrt{n-m_n}}\right)=\frac{1}{m_n}\sum_{\ell=1}^{K}X_{\ell}^2(\omega)\ind_{\{\vert X_{\ell}(\omega)\vert> \sqrt{\ell}\}}\left(1\wedge \frac{\vert X_{\ell}(\omega)\vert}{\sqrt{n-m_n}}\right)\converge{9}{n}{+\infty}{}0.
\end{equation}
In the other part, we have 
\begin{align*}
	\mathbb E\!\left[\sum_{\ell=1}^{+\infty}\vert X_{\ell}\vert^3 \ell^{-3/2}\ind_{\vert X_{\ell}\vert\leqslant \sqrt{\ell}}\right]
	&=\sum_{\ell=1}^{+\infty}\mathbb E\!\left[\vert \xi_1\vert^3\ell^{-3/2}\ind_{\vert \xi_1\vert\leqslant\sqrt{\ell}}\right]
	=\mathbb E\!\left[\vert\xi_1\vert^3\sum_{\ell\geqslant\xi_1^2}\ell^{-3/2}\right]\lesssim\mathbb E[\xi_1^2]<+\infty.
\end{align*}
Consequently, $\sum_{\ell=1}^{+\infty}\vert X_{\ell}\vert^3 \ell^{-3/2}\ind_{\vert X_{\ell}\vert\leqslant \sqrt{\ell}}<+\infty$ a.s. So, using Kronecker's lemma, we derive 
\begin{equation}\label{control_A_1_n_part2}
\frac{1}{m_n}\sum_{\ell=1}^{m_n}X_{\ell}^2\ind_{\{\vert X_{\ell}\vert\leqslant  \sqrt{\ell}\}}\left(1\wedge \frac{\vert X_{\ell}\vert}{\sqrt{n-m_n}}\right)\leqslant \frac{\sum_{\ell=1}^{m_n}\vert X_{\ell}\vert^3\ind_{\vert X_\ell\vert\leqslant\sqrt{\ell}}}{m_n\sqrt{n-m_n}}\lesssim\frac{\sqrt{\theta}\sum_{\ell=1}^{m_n}\vert X_{\ell}\vert^3\ind_{\vert X_\ell\vert\leqslant\sqrt{\ell}}}{m_n^{3/2}\sqrt{1-\theta}}\converge{14}{n}{+\infty}{a.s.}0.
\end{equation}
Combining \eqref{control_A_1_n_part1} and \eqref{control_A_1_n_part2}, we get 
$$
\frac{1}{m_n}\sum_{\ell=1}^{m_n}X_{\ell}^2\left(1\wedge \frac{\vert X_{\ell}\vert}{\sqrt{n-m_n}}\right)\converge{14}{n}{+\infty}{a.s.}0.
$$
Consequently $A_{1,n}\converge{12}{n}{+\infty}{a.s.}0$ and $(n-m_n)\mathbb E_{\F_{m_n}}[\vert R_n\vert]\converge{12}{n}{+\infty}{a.s.}0$. Keeping in mind that 
$$
\mathbb E_{\F_{m_n}}[\exp\left(it\sqrt{\kappa_n}H_n\right)]=\left(1-\frac{t^2\kappa_n\mathbb E_{\F_{m_n}}[Y'^2_{m_n+1}]}{2(n-m_n)}+\mathbb E_{\F_{m_n}}[R_n]\right)^{n-m_n},
$$
we get
$$
\mathbb E_{\F_{m_n}}[\exp\left(it\sqrt{\kappa_n}H_n\right)]\converge{19}{n}{+\infty}{$\mathbb L^1$ and a.s.}\exp\left(-\frac{\theta(1-\theta)\sigma^2t^2}{2}\right).
$$
Consequently, the convergence \eqref{Lemme_clef_convergence_fonction_caracteristique_conditionnelle} holds for any $-1\leqslant a<1$.\\
\\
Now, we assume that $-1\leqslant a<1/2$. Using Theorem \ref{TLC}, we have
\begin{equation}\label{conv_fonction_caracteristique_G_n}
\lim_{n\to+\infty}\mathbb E[\exp(itG_n)]=\exp\left(-\frac{\tau^2\sigma^2t^2}{2(1-2a)}\right)\quad\textrm{where}\quad\tau=\lim_{n\to+\infty}\tau_n=\theta+a(1-\theta).
\end{equation}
Combining Lemma \ref{Lemme_clef_convergence_fonction_caracteristique_conditionnelle} and \eqref{conv_fonction_caracteristique_G_n}, we obtain $\lim_{n\to+\infty}\psi_n(t)=\exp\left(-\frac{\sigma^2t^2}{2}\left(\frac{\tau^2}{1-2a}+\theta(1-\theta)\right)\right)$. So, for $-1\leqslant a<1/2$, we have 
$$
\sqrt{m_n}\left(\frac{ S_n}{n}-\frac{(1-\alpha)\mu_1}{1-a}\right)\converge{15}{n}{+\infty}{\textrm{Law}}\mathcal N\left(0,\frac{\tau^2\sigma^2}{1-2a}+\theta(1-\theta)\sigma^2\right).
$$
Assume that $a=1/2$. If we denote $G'_n:=(\log m_n)^{-1/2} G_n$ and $H'_n:=(\log m_n)^{-1/2} H_n$ where $G_n$ and $H_n$ are defined by \eqref{def_G_n_et_H_n}, then 
$$
\frac{\sqrt{m_n}}{\sqrt{\log m_n}}\left(\frac{ S_n}{n}-2(1-\alpha)\mu_1\right)=G'_n+\sqrt{\kappa_n}H'_n
$$
with 
\begin{align*}
\mathbb E\left[\kappa_nH'^2_n\right]&=\frac{\kappa_n}{(n-m_n)\log m_n}\mathbb E\left[\left(\sum_{k=m_n+1}^nY'_k\right)^2\right]\\
&=\frac{\kappa_n}{(n-m_n)\log m_n}\mathbb E\left[\sum_{k=m_n+1}^n\mathbb E_{\F_{m_n}}[Y'^2_k]+2\sum_{m_n<k<\ell\leqslant n}\underbrace{\mathbb E_{\F_{m_n}}[Y'_kY'_{\ell}]}_{=0}\right]\\
&=\frac{\kappa_n\mathbb E[Y'^2_{m_n+1}]}{\log m_n}\converge{15}{n}{+\infty}{}0.
\end{align*}
Moreover, since $a=1/2$, using Theorem \ref{TLC}, we have 
$$
G'_n\converge{15}{n}{+\infty}{\textrm{Law}}\mathcal N\left(0,\mu_2\tau^2-4\tau^2(1-\alpha)^2\mu_1^2\right)\quad\textrm{with}\quad\tau=\theta+a(1-\theta).
$$
Consequently, using Slutsky's lemma, we obtain 
$$
\frac{\sqrt{m_n}}{\sqrt{\log m_n}}\left(\frac{ S_n}{n}-2(1-\alpha)\mu_1\right)\converge{19}{n}{+\infty}{\textrm{Law}}\mathcal N\left(0,\mu_2\tau^2-4\mu_1^2\tau^2(1-\alpha)^2\right).
$$
Assume that $1/2<a<1$. As before, one can notice that 
$$
m_n^{1-a}\left(\frac{S_n}{n}-\frac{(1-\alpha)\mu_1}{1-a}\right)=G''_n+\sqrt{\kappa_n}H''_n\quad\textrm{with}\quad G''_n:=m_n^{\frac{1}{2}-a}G_n\quad\textrm{and}\quad H''_n:=m_n^{\frac{1}{2}-a}H_n.
$$
Using the martingale convergence theorem (see \cite{Duflo1997}), we have $$
M_{m_n}:=a_{m_n}\tilde{T}_{m_n}=\frac{(\alpha-a)\mu_1}{1-a}+\sum_{\ell=1}^{m_n}a_{\ell}\varepsilon_{\ell}\converge{12}{n}{+\infty}{a.s.}M_{\infty}:=\frac{(\alpha-a)\mu_1}{1-a}+\sum_{\ell=1}^{+\infty}a_{\ell}\varepsilon_{\ell}.
$$ 
So, if we denote $L:=\frac{M_{\infty}}{\Gamma(a+1)}$, then 
$$
\sqrt{m_n^{2a-1}}\left(m_n^{1-a}\left(\frac{S_n}{n}-\frac{(1-\alpha)\mu_1}{1-a}\right)-\tau_nL\right)=G_n-\tau_n\sqrt{m_n^{2a-1}}L+\sqrt{\kappa_n}H_n=J_n+\sqrt{\kappa_n}H_n
$$
where 
$$
J_n:=
\tau_n\sqrt{m_n^{2a-1}}\left(m_n^{1-a}\left(\frac{T_{m_n}}{m_n}-\frac{(1-\alpha)\mu_1}{1-a}\right)-L\right).
$$
Let $t\in\R$ be fixed and define $\tilde{\psi}_n(t):=\mathbb E[\exp(it\left(J_n+\sqrt{\kappa_n}H_n\right)]=\mathbb E[\exp(itJ_n)\mathbb E_{\F_{m_n}}[\exp(it\sqrt{\kappa_n}H_n)]]$. Using Theorem \ref{TLC} and Lemma \ref{Lemme_clef_convergence_fonction_caracteristique_conditionnelle}, we have 
$$
J_n\converge{15}{n}{+\infty}{\textrm{Law}}\mathcal N\left(0,\frac{\tau^2\sigma^2}{2a-1}\right)
\quad\textrm{and}\quad
\mathbb E_{\F_{m_n}}[\exp\left(it\sqrt{\kappa_n}H_n\right)]\converge{18}{n}{+\infty}{$\mathbb L^1$ and a.s.}\exp\left(-\frac{\theta(1-\theta)\sigma^2t^2}{2}\right). 
$$
So, for any $t\in\R$, 
$$
\lim_{n\to+\infty}\tilde{\psi}_n(t)=\exp\left(-\frac{\sigma^2t^2}{2}\left(\frac{\tau^2}{2a-1}+\theta(1-\theta)\right)\right).
$$
Finally, if $1/2<a<1$, then 
$$
\sqrt{m_n^{2a-1}}\left(m_n^{1-a}\left(\frac{S_n}{n}-\frac{(1-\alpha)\mu_1}{1-a}\right)-\tau_nL\right)\converge{19}{n}{+\infty}{\textrm{Law}}\mathcal N\left(0,\frac{\tau^2\sigma^2}{2a-1}+\theta(1-\theta)\sigma^2\right).
$$
The Proof of Theorem \ref{TLC_bis} is complete.$\hfill\Box$\\
\\
\textbf{Acknowledgements}.  The authors would like to express their sincere gratitude to the anonymous referee for their careful reading of the manuscript and for their valuable comments and suggestions. In particular, we are deeply grateful for the referee's insightful recommendation that led to a substantial simplification of the proof of Theorem $3$ $\textrm{iii)}$.

\bibliographystyle{plain}
\bibliographystyle{stylename}
\bibliography{/Users/elmachko/Documents/Documents_de_travail/Enseignement_et_Recherche/Recherche/BIB/xbib.bib}
\end{document}